\documentclass{amsart}
\usepackage{amssymb,amsfonts}
\usepackage[all,arc]{xy}
\usepackage{enumerate}
\usepackage{mathrsfs}
\usepackage{graphicx}
\usepackage{mathtools}
\usepackage{leftidx}
\usepackage{bm}
\newtheorem{thm}{Theorem}[section]
\newtheorem{cor}[thm]{Corollary}
\newtheorem{prop}[thm]{Proposition}
\newtheorem{lem}[thm]{Lemma}

\theoremstyle{definition}
\newtheorem{defn}[thm]{Definition}

\newtheorem{exmp}[thm]{Example}

\theoremstyle{remark}
\newtheorem{rem}[thm]{Remark}

\makeatletter
\let\c@equation\c@thm
\makeatother
\numberwithin{equation}{section}

\def\bthm{\begin{thm}}
\def\ethm{\end{thm}}
\def\blm{\begin{lem}}
\def\elm{\end{lem}}
\def\bdf{\begin{defn}}
\def\edf{\end{defn}}
\def\bpf{\begin{proof}}
\def\epf{\end{proof}}
\def\bpp{\begin{prop}}
\def\epp{\end{prop}}
\def\bcor{\begin{cor}}
\def\ecor{\end{cor}}
\def\brm{\begin{rem}}
\def\erm{\end{rem}}
\def\D{{D}}
\def\A{{\mathscr A}}
\def\R{{\mathbb R}}
\def\F{{\mathbb F}}
\def\L{\boldsymbol{L}}

\def\bZ{\mathbb{Z}}

\def\bR{{\mathbb R}}

\def\bF{\mathbb{F}}
\def\cA{\mathscr{A}}

\def\cC{\mathcal{C}}

\def\cS{\mathcal{S}}

\def\cL{\mathcal{L}}
\def\cQ{\mathcal{Q}}

\newcommand{\Alg}{\mathsf{Alg}}
\newcommand{\dgcat}{\mathsf{dgCat}}

\newcommand{\Cat}{\mathsf{Cat}}

\newcommand{\Set}{\mathsf{Set}}
\newcommand{\HoC}{{\rm Ho}(\mathcal{C})}
\newcommand{\Ho}{\textrm{Ho}}

\newcommand{\Mod}{\textrm{Mod}}

\newcommand{\Gr}{\mathsf{Gr}}

\newcommand{\Birack}{\mathsf{Birack}}
\newcommand{\Biquandle}{\mathsf{Biquandle}}
\newcommand{\Perv}{\mathsf{Perv}}

\newcommand{\hatC}{\widehat{\mathcal{C}}}
\newcommand{\hatA}{\hat{A}}
\newcommand{\Ob}{\rm{Ob}}
\newcommand{\Color}{\textrm{Col}}

\newcommand{\Ch}{\rm{Ch}}
\newcommand{\sGr}{\mathsf{sGr}}
\newcommand{\Fuk}{\mathsf{Fuk}}
\newcommand{\Sh}{\mathsf{Sh}}
\newcommand{\Aug}{\mathsf{Aug}}

\newcommand{\Hom}{{\rm Hom}}

\newcommand{\Aut}{{\rm Aut}}

\newcommand{\id}{{\rm id}}
\newcommand{\colim}{{\rm colim}}

\newcommand{\wa}{{\rm wa}}

\newcommand{\cyl}{{\rm Cyl}}
\newcommand{\coker}{{\rm coker}}

\newcommand{\Tw}{{\rm Tw}}

\newcommand{\xhocl}[2]{
\bm{L}\rm{coeq}
\bigl[\xymatrix{
#1
\ar@<0.5ex>[r]^{#2}   \ar@<-0.5ex>[r]_{\id}
& #1
} \bigr]
}

\newcommand{\xcatcl}[2]{
{\rm coeq}
\bigl[ \xymatrix{
#1
\ar@<0.5ex>[r]^{#2}   \ar@<-0.5ex>[r]_{\id}
& #1
} \bigr]
}

\newcommand{\xcoeq}[4]{
{\rm coeq}
\bigl[\xymatrix{
#1
\ar@<0.5ex>[r]^{#2}   \ar@<-0.5ex>[r]_{#3}
& #4
} \bigr]
}

\newcommand{\xhocoeq}[4]{
{\rm hocoeq}
\bigl[ \xymatrix{
#1
\ar@<0.5ex>[r]^{#2}   \ar@<-0.5ex>[r]_{#3}
& #4
} \bigr]
}

\newcommand{\into}{\,\hookrightarrow\,}

\newcommand{\onto}{\,\twoheadrightarrow\,}

\newcommand{\sonto}{\,\stackrel{\sim}{\twoheadrightarrow}\,}
\newcommand{\ra}{\,\rightarrow \,}
\newcommand{\la}{\,\leftarrow \,}

\newcommand{\YB}{\sigma : A \amalg A \ra A \amalg A}

\newcommand{\xra}[1]{\, \xrightarrow{#1} \,}
\newcommand{\xla}[1]{\, \xleftarrow{#1} \,}
\newcommand{\xinto}[1]{\,\stackrel{#1}{\hookrightarrow}\,}

\newcommand{\infinity}{1}

\DeclareMathAlphabet{\mathpzc}{OT1}{pzc}{m}{it}

\title{Perverse sheaves and knot contact homology}
\author{Yuri Berest}
\address{Department of Mathematics,
Cornell University, Ithaca, NY 14853-4201, USA}
\email{berest@math.cornell.edu}
\author{Alimjon Eshmatov}
\address{Department of 
Mathematics and Statistics,
University of Toledo,
Toledo, OH 43606-3390, USA}
\email{alimjon.eshmatov@utoledo.edu}
\author{Wai-Kit Yeung}
\address{Department of Mathematics,
Cornell University, Ithaca, NY 14853-4201, USA}
\email{wy236@cornell.edu}
\begin{document}

\begin{abstract}
In this paper, which is a research announcement of results that will appear in \cite{BEY},
we give a new algebraic construction of knot contact homology 
in the sense of Ng \cite{Ng05a}. For a link $L$ in $ {\mathbb R}^3 $, we define
a differential graded (DG) $k$-category $ \tilde{\A}_L $ with finitely many objects, whose quasi-equivalence class is a topological invariant of $ L $. In the case when $L$ is a knot, the endomorphism algebra of a distinguished object of $ \tilde{\A}_L $ coincides with the  fully noncommutative knot DGA  as defined by Ekholm, Etnyre, Ng and Sullivan in \cite{EENS13a}. The input of our construction is a natural action of the braid group $B_n$ on the category of perverse sheaves on a two-dimensional disk with singularities at $n$ marked points, studied by Gelfand, MacPherson and Vilonen in \cite{GMV96}. As an application, we show that the category of finite-dimensional representations of the link $k$-category $ \tilde{A}_L = H_0(\tilde{\A}_L) $ defined as the $0$-th homology of  $ \tilde{\A}_L $ is equivalent to the category of perverse sheaves on $ \bR^3 $ which are singular along the link $ L $.  We also obtain several generalizations of the category  $ \tilde{\A}_L $ by extending the Gelfand-MacPherson-Vilonen braid group action.
\end{abstract}

\maketitle

\tableofcontents

\section{Introduction}

In a series of papers \cite{Ng05a, Ng05b, Ng08, Ng11, Ng14}, L.~Ng introduced and studied a new algebraic invariant of a link $L$ in $ \bR^3 $ represented by
a semi-free differential graded (DG) algebra $ {\mathcal A}_L $. The structure of this DG algebra (referred to {\it op.cit.} as a combinatorial knot DGA) is determined by an element of a braid group $B_n$ representing the link $L$. The homology of 
$ {\mathcal A}_L$ is called the {\it knot contact homology} 
$ HC_\ast(L) $ as it coincides with the Legendrian contact homology\footnote{in the sense of \cite{Eli98} 
(see also \cite{EES05, EES07})}
of the unit conormal bundle $ \Lambda_L \subseteq ST^*\bR^3 $ of $L$.
In fact, it was conjectured in \cite{Ng05a, Ng05b} and recently proved in \cite{EENS13a, EENS13b} that the entire combinatorial knot DGA is isomorphic to a geometrically defined DG algebra representing
the Legendrian contact homology of $ \Lambda_L $. 

Our original motivation  was to understand Ng's combinatorial proof of the invariance of ${\mathcal A}_L $ (up to stable isomorphism) under the Markov moves. We should remark that, although the differential of ${\mathcal A}_L $ is defined in \cite{Ng05a} by an explicit formula, its combinatorial structure is fairly complicated and its algebraic origin seems mysterious. Even the fact that the $0$-th homology of $ {\mathcal A}_L $ is a link invariant is far from being obvious from the definition of 
\cite{Ng05a} ({\it cf.} \cite[Section~4.3]{Ng05a}).
As a result, we have come up with a different, more conceptual construction
that makes the Markov invariance of $ {\mathcal A}_L $ quite transparent\footnote{In fact, the invariance of our construction under type
I Markov moves follows directly from its definition.} 
and, more importantly, places knot contact homology in one row with other classical invariants, such as knot groups and Alexander modules.

To clarify the ideas we begin by recalling a classical theorem of E. Artin and J. Birman that 
gives a natural presentation of the link group $ \pi_1(\R^3\!\setminus\! L) $ in terms of a braid representing $L$. 
Let $ \D $ be the unit disk in $\R^2$, and let $ \{p_1, \ldots, p_n\} \subset \D $ be a set of
distinct points in the interior of $D$. It is well known that the braid group on $n$-strands,  $ B_n $, can be identified with the mapping class group of $ (\D\!\setminus\! \{p_1, \ldots, p_n\}, \partial \D) $, and as such it acts 
naturally on the fundamental group $ \pi_1(\D\!\setminus\! \{p_1, \ldots, p_n\}, p_0) $, where 
$ p_0 \in \D\!\setminus \!\{p_1, \ldots, p_n\} $ is a basepoint which we choose near
the boundary of $ \D $. The fundamental group 
$ \pi_1(\D\!\setminus\! \{p_1, \ldots, p_n\}, p_0) $ is
a free group $ \F_n $ of rank $n$ based on generators $ x_1, \ldots, x_n $ which correspond to
small loops in $ \D\!\setminus\! \{p_1, \ldots, p_n\} $ around the points $ p_i $. 
Explicitly, in terms of these generators, the action of $B_n$ on
$ \pi_1(\D\!\setminus\! \{p_1, \ldots, p_n\}, p_0) \cong \F_n $ is given by
\begin{equation}\label{artinact}
\sigma_i\,: \left\{
\begin{array}{lll}
x_i & \mapsto & x_i \,x_{i+1}\, x_i^{-1} \\
x_{i+1} & \mapsto  & x_i\\
x_j & \mapsto & x_j \quad (j \not= i, i+1) 
\end{array}\right.
\end{equation}
where $ \sigma_i $ $ (i=1, 2, \ldots, n-1) $ are the standard generators of $B_n$. This action is usually
called the {\it Artin representation}  as it provides a faithful realization of $B_n$ as a subgroup of
$ \Aut(\F_n) $. Now, the Artin-Birman Theorem (see \cite[Theorem~2.2]{Bir74})
asserts that the fundamental group of the complement of the link $ L = \hat{\beta} \subset \R^3 $
corresponding to a braid $ \beta \in B_n $ has the presentation
\begin{equation}
\label{artpres}
\pi_1(\R^3\!\setminus\!L)\, \cong\, \langle x_1, \,x_2,\, \ldots\,,\, x_n\ |\ \beta(x_1) = x_1,\, \beta(x_2) = x_2,\, 
\ldots,\, \beta(x_n) = x_n \rangle\ ,
\end{equation}
where $ \beta(x_i) $ denotes the action of $ \beta $ on $ x_i$ via the Artin representation.

We abstract this situation in the following way. Let $\cC $ be a category with finite colimits.
We assume that we are given a family of braid group actions $\,\varrho_n: B_n \to \Aut\,A^{(n)}\,$,
$\,n \ge 1 \,$, on objects of $ \cC $ having the properties:
\begin{enumerate}
\item[(1)] For each $ n \ge 1 $, $ \, A^{(n)} $ is the $n$-fold coproduct of 
one and the same object $ A $ of $\cC$.
\item[(2)] The actions $\varrho_n $ are \emph{local} and {\it homogeneous} 
in the sense that each $\sigma_i \in B_n $ acts only on the $(i,i+1)$-copy of $A^{(2)} $ in $A^{(n)}$ 
while keeping the rest fixed, and any two standard generators of $ B_n $ 
act in the same way on the corresponding copies of $ A^{(2)} $ for all $n\ge 1$.
\end{enumerate}
Such braid group actions are determined (generated) by a single morphism $ \sigma: A \amalg A \to A\amalg A $ in the category
$ \cC $ that we call a \emph{cocartesian Yang-Baxter operator} ({\it cf.} Definition~\ref{cocartesian_Yang_Baxter_operators_definition}). For example, the Artin representations are generated by
a cocartesian Yang-Baxter operator in the category of groups given by 
$ \sigma: \F_1 \amalg \F_1 \to \F_1 \amalg \F_1 $, $\,x_1 \mapsto x_1 x_2 x_1^{-1}$,$\  x_2 \mapsto x_1 $.

Now, for an arbitrary cocartesian Yang-Baxter operator $ (A,\sigma) $, we define a universal construction $ {\mathcal L}(A, \sigma) $ that associates to each braid $\beta \in B_n $ the coequalizer of the endomorphisms $ \id $ and $ \beta $ of the object $ A^{(n)} $, or equivalently, the following pushout in $ \cC \,$:
\begin{equation}
\label{catclos}
{\mathcal L}(A, \sigma)[\beta] := \xcatcl{\!\!A^{(n)}\!\!}{\beta} \,
= \,{{\rm colim}}\,[\,A^{(n)} \xleftarrow{\,(\beta,\, \id)\,} A^{(n)} \amalg A^{(n)} \xrightarrow{\,(\id,\, \id)\,} A^{(n)} \,] \ .
\end{equation}
We call \eqref{catclos} the \emph{categorical closure of the braid $ \beta $ on the object} $A$ with respect to the Yang-Baxter operator $ \sigma $. This terminology can be justified by the following ``picture''\footnote{In fact, this ``picture'' of a categorical braid closure can be formalized by using the diagrammatic tensor calculus developed by A. Joyal, R.~Street and others (see \cite{FY92, Shu94, JS91a, JS91b, JS93, JSV96, RT90, Tur10}). We briefly discuss it at the end of Section \ref{Yang_Baxter_operators_section}.} of the pushout \eqref{catclos} that
manifestly exhibits it as ``a braid closure on $A$'':
\[ 
\includegraphics[scale=0.3]{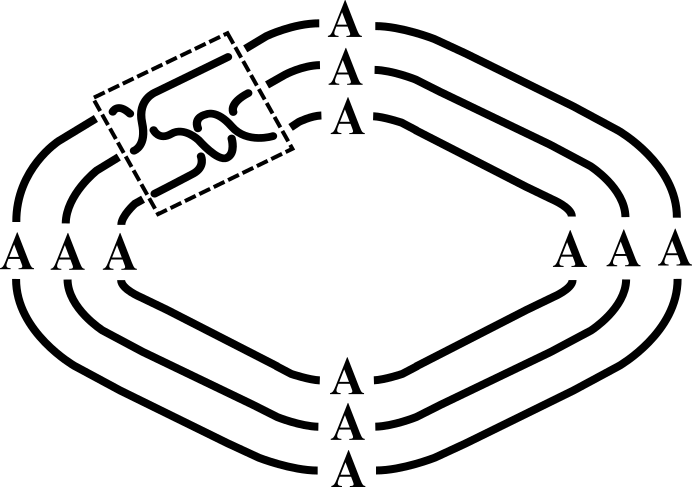}
\]
\[ 
\colim\,[\,A^{(3)} \xleftarrow{(\beta,\, \id)} A^{(3)} \amalg A^{(3)} \xrightarrow{(\id, \,\id)} A^{(3)} \,] \qquad \,
\]

In the case of Artin representations, the Artin-Birman Theorem \eqref{artpres}  implies that
$\, {\mathcal L}(\F_1, \sigma)[\beta] \cong \pi_1(\R^3\!\setminus\!L)\,$. This means, in particular, that 
$  {\mathcal L}(\F_1, \sigma)[\beta] $ is a link invariant.

In general, we show that, if a cocartesian Yang-Baxter operator $(A,\sigma)$ satisfies some natural conditions, which
we call the {\it Reidemeister conditions} (see Defintion~\ref{Reidemeister_operators_definition}),
then the isomorphism class of the categorical closure of any braid with respect to $ (A,\sigma) $
is stable under the Markov moves, and hence defines a link invariant ({\it cf.} Theorem \ref{Wada_theorem_trivial_torsion} and Theorem \ref{Wada_theorem_nontrivial_torsion}). Apart from the group $ \pi_1(\R^3\!\setminus L) $, many classical link
invariants arise in this way (see, for example, Theorem~\ref{Alexm} that represents as a categorical braid 
closure the Alexander module).

Next, we consider the category $\Perv(D, \{p_1,\ldots,p_n\})$ of perverse sheaves on the disk $ D $ with only possible singularities at the points $\{p_1, \ldots , p_n\} $. In \cite{GMV96}, Gelfand, MacPherson and Vilonen showed 
that $\Perv(D, \{p_1,\ldots,p_n\})$ is equivalent to the category $  \tilde{\mathscr Q}^{(n)} $ of finite-dimensional 
$k$-linear representations of the following quiver
%
%
%
\begin{equation*}
Q^{(n)} = \qquad
\vcenter{
\xymatrix@C=2em@R=1.5em{
& & 1 \ar@/^0.5pc/[dd]|{a_1} \\
n \ar@/^0.5pc/[rrd]|{a_n}  & & & & 2 \ar@/^0.5pc/[lld]|{a_2} \\
 & & 0 \ar@/^0.5pc/[uu]|{a_1^*} \ar@/^0.5pc/[rru]|{a_2^*} \ar@/^0.5pc/[rrd]|{a_3^*} \ar@/^0.5pc/[llu]|{a_n^*} \ar@{}[dl]|{\ddots} \\
&  &  & & 3 \ar@/^0.5pc/[llu]|{a_3} 
}
}
\end{equation*}
such that the operators $ T_i := e_0 + a_i a_i^* $ act as isomorphisms
for all $ i = 1,\,2,\,\ldots, n\, $. More formally,  $ \tilde{\mathscr Q}^{(n)} $
can be described as the category $ \Mod\,	\tilde{A}^{(n)} $
of finite-dimensional modules over the $k$-category
\begin{equation}  
\label{widetilde_A}
	\tilde{A}^{(n)} := k \langle Q^{(n)} \rangle
[T_1^{-1},\,\ldots\,,\,T_n^{-1}]
\end{equation} 
which is obtained by localizing the path category of $ Q^{(n)} $ at 
the set of morphisms $\, \{T_1, \ldots, T_n\} \,$.
Now, the braid group $B_n$ acts on the disk $ D $ with $n$ marked points $ \{p_1, \ldots ,p_n\} $ as a mapping class group, and this naturally induces an action on the category  $\Perv(D, \{p_1, \ldots ,p_n\}) $.
It was shown in \cite{GMV96} that, under the equivalence   $\,\Perv(D, \{p_1,\ldots,p_n\}) \simeq  \tilde{\mathscr Q}^{(n)} $, the action of $B_n$ on the category of perverse sheaves corresponds to a \emph{strict} action on the category $\tilde{\mathscr{Q}}_n $ ({\it cf.} \cite[Proposition~1.3]{GMV96}). This, in turn, induces an action of $B_n$ on the $k$-category $	\tilde{A}^{(n)} $, which is given explicitly (on generating morphisms of $	\tilde{A}^{(n)} $) by the following formulas
\begin{equation} 
\label{GMV}
\sigma_i \, : \,  \begin{cases} 
     a_{i} \, \mapsto \, T_{i} \,a_{i+1} &  \\
     a_{i+1} \, \mapsto \, a_{i} & \\ 
     a_j \, \mapsto \, a_{j} &  (j \neq i, i+1) \\
     a_{i}^* \, \mapsto \,  a_{i+1}^* \, T_{i}^{-1}&  \\
     a_{i+1}^* \, \mapsto \, a_{i}^* & \\
      a_j^* \, \mapsto \, a_j^* & (j \neq i, i+1) 
   \end{cases}  
\end{equation}
We call \eqref{GMV} the {\it Gelfand-MacPherson-Vilonen (GMV) braid action}. 

The GMV braid actions are generated by a single cocartesian Yang-Baxter operator in the category of (small) {\it pointed} $k$-categories  $ \Cat_k^* $. Specifically, for each $ n\ge 1 $, the $k$-category $    	\tilde{A}^{(n)} $ is the coproduct (fusion product) in $ \Cat_k^* $ of $ n $ copies of the $k$-category $ \tilde{A} = k\langle Q \rangle [T^{-1}] $, where $ k\langle Q \rangle $ is the path category of the quiver
$\, Q = \bigl[  
\vcenter{
\xymatrix@1{
\infinity  \ar@/^/[rr]^{a}  
&& 0  \ar@/^/[ll]^{a^{*}} 
}
} \bigr]
\,$
with the distinguished object $0$. The corresponding Yang-Baxter operator
$\,\sigma: \tilde{A} \amalg  \tilde{A} \to \tilde{A} \amalg  \tilde{A} $ is given by
\begin{equation} 
\label{GMV1}
(a_1,\,a_1^*) \,\mapsto\, (T_1 \,a_2,\, a_2^*\, T_1^{-1})\ ,\quad (a_2,\,a_2^*) \,\mapsto\, (a_1,\,a_1^*)\ .
\end{equation}
Just as in the case of Artin actions, it is easy to check that \eqref{GMV1} satisfies the Reidemeister conditions, and 
hence the categorical braid closure with respect to $ (\tilde{A}, \sigma) $ is a link invariant. For a given 
$ \beta \in B_n $, this invariant is represented by the equivalence class of the $k$-category 
$ \tilde{A}_L := {\mathcal L}(\tilde{A}, \sigma)[\beta] $, which we call the {\it (fully noncommutative) 
link $k$-category}\footnote{Strictly speaking, the categorical braid closure gives a specialization of the fully noncommutative link $k$-category, with all longitude parameters set to be $1$ (see Remark~\ref{catcl_of_GMV_is_not_cord_cat}). For a general definition of $ \tilde{A}_L$, we refer to Section~\ref{FN}, Definition~\ref{FNCDG}.} of $ L = \hat{\beta} $. In Section~\ref{FN}, we will show that the 
$k$-category $ \tilde{A}_L $ is a natural extension of the fully noncommutative  cord algebra of \cite{EENS13a, Ng14} 
in the sense that the latter can be identified with the endomorphism algebra of an object in $ \tilde{A}_L $. 
Thus, in our algebraic formalism, 
the link category $ \tilde{A}_L $ arises exactly the same way as the link group $ \pi_1(\R^3\!\setminus\!L)$,
provided we take as an input the Gelfand-MacPherson-Vilonen braid action instead of
the Artin representation.

At this point, we pause to remark that the notion of a categorical braid closure has already appeared in the literature: explicitly -- in the case of groups (see \cite{Wad92, CP05}), and  in a somewhat disguised form, in the theory of quandles 
(see, for example, \cite{FJK04, CSWES09}). From this last perspective, our results give a precise interpretation of 
such geometric knot invariants as a cord algebra in combinatorial terms of racks and quandles (see Remark~\ref{relation_with_quandles_remark} below).

However, our main observation is that the simple categorical formalism we outlined above admits an interesting generalization to homotopical contexts. Specifically, if the category $\cC$ that we work with has 
a natural class $ {\mathcal W} $ of weak equivalences (e.g., $\cC$ is a Quillen model category or a homotopical category in the sense of \cite{DHKS04}), then the operation of a categorical braid closure is usually not invariant under weak 
equivalences, {\it i.e.} it is not well-defined\footnote{The problem is that pointwise weak equivalences of diagrams do not necessarily induce weak equivalences of colimits, so the objects defined by colimits of diagrams defined up to homotopy are not well-defined, even up to homotopy type.} in the homotopy category $ \HoC = \cC[{\mathcal W}^{-1}] $.
In abstract homotopy theory, there is a standard way to remedy this problem: namely, replace a homotopy non-invariant
functor $F$ by its derived functor $ \L F $, which gives a universal approximation to $F$ on the level of homotopy
categories (see, e.g., \cite[Section~9]{DS95}). In our situation, we can define a ``derived'' version of the categorical braid closure by simply replacing the `$ \colim $' in the definition \eqref{catclos} by its derived functor: 
the homotopy colimit `$ \L\colim $'.
To be precise, given a cocartesian Yang-Baxter operator $ (A, \sigma) $ in (say) a model category $ \cC $, we define the {\it homotopy braid closure} of  $ \beta \in B_n $ with respect to $ (A, \sigma) $ by
\begin{equation}
\label{hobrcl}
\L {\mathcal L}(A, \sigma)[\beta] := 
\L{\colim} [A^{(n)} \xleftarrow{(\beta,\,\id)} A^{(n)}\amalg A^{(n)} 
\xrightarrow{(\id,\,\id)} A^{(n)}] \ .
\end{equation}
One of our main results (Theorem~\ref{invHBC}) states that if $ (A, \sigma) $ satisfies the Reidemeister conditions (and $A$ is flat in an appropriate sense), then the weak equivalence class of the homotopy braid closure on $A$, {\it i.e.} the isomorphism class of \eqref{hobrcl} in the homotopy category $ \HoC $, is invariant under the Markov moves, and hence defines a link invariant.
This last invariant is more refined than the one given by the usual categorical braid closure in the same way as 
the homotopy type of a topological space is a more refined invariant of the space than just its fundamental group.

Now, let us return to our basic example of the cocartesian Yang-Baxter operator 
$ (\tilde{A}, \sigma) $ associated to the GMV action, see \eqref{GMV1}. 
To define the homotopy braid closure with respect to this operator we will regard the $k$-category $ \tilde{A}$ 
as an object of the category $\dgcat^*_k$ comprising all (small) pointed DG categories. 
The category $\dgcat^*_k$ has a natural model structure, in which the weak equivalences are
the quasi-equivalences\footnote{Recall that a {\it quasi-equivalence} of DG categories is 
a DG functor $F: {\mathscr A} \to {\mathscr B} $ such that $ F:\, {\mathscr A}(X,Y) \to {\mathscr B}(FX, FY) $
is a quasi-isomorphism of $k$-complexes for all objects $X,Y \in {\rm Ob}({\mathscr A}) $ and the induced functor
on the $0$th homology  $ H_0(F):\, H_0({\mathscr A}) \stackrel{\sim}{\to} H_0({\mathscr B}) $ is an 
equivalence of categories.} of DG categories (see \cite{Tab05}). Its homotopy theory has been  
extensively studied in recent years with a view towards applications in algebraic geometry 
and algebraic $K$-theory (see, e.g., \cite{Ke06, Toen07, Toen11} and references therein).

The homotopy braid closure of the GMV action in the model category $\dgcat^*_k$ gives a new link invariant,
which is a quasi-equivalence class of DG categories. For a given $ \beta \in B_n $, formula \eqref{hobrcl}
allows us, in fact, to construct an explicit representative for the corresponding quasi-equivalence class 
that we call the \emph{fully noncommutative link DG category} $\,\tilde{\cA}_L $ (see Definition~\ref{FNCDG}).
If we assume, for simplicity, that $ L $ is a knot ({\it i.e.}, a link with a single component), then
$\tilde{\cA}_L $ contains a distinguished object, and the endomorphism DG algebra of that object is isomorphic 
to the fully noncommutative knot DGA constructed in \cite{EENS13a}. This observation is part of 
Theorem~\ref{fully_noncommutative_link_DG_category_main_theorem} that we state in full generality (for 
links with an arbitrary number of components) but do not prove in this paper. Instead,
we sketch a proof of an analogous result -- Theorem~\ref{endomorphism_DGA_at_1_recovers_knot_DGA_theorem} --
that identifies the {\it framed} knot DGA  (originally introduced in \cite{Ng08}) with
the DG endomorphism algebra of a distinguished object in the homotopy braid closure of a
modified GMV action. The modification amounts to collapsing all objects of the GMV $k$-category
$ \tilde{A}^{(n)} $, except for the base object `$0$', to a single object `$1$',
while preserving all the generating morphisms $ a_i $ and $ a_i^* $. 
We also impose $n$ extra relations $\,a_i^*\,a_i = (\mu - 1)\,e_1\,$, one for each $i=1, 2, \ldots, n$,
that depend on an invertible central parameter $ \mu $ in the ground ring $k$.
The resulting $k$-category $ A^{(n)} $ with two objects $ \{0, 1\} $ inherits the GMV braid action 
\eqref{GMV}, 
and one can still define its homotopy braid closure by formula \eqref{hobrcl}\footnote{To introduce
the second central parameter $ \lambda \in k^{\times} $ we also modify the arrow $ (\beta, \id) $ 
in the homotopy colimit \eqref{hobrcl} by appropriately twisting the action map 
$ \beta: A^{(n)} \to A^{(n)} $ (see Section \ref{The_GMV_operator_section}).}.

Now, as in the case of topological spaces, to compute the homotopy colimit of a diagram like
\eqref{hobrcl} one should first `resolve' the objects by their cofibrant models
$ R \sonto A^{(n)} $, then replace one of the arrows by a (weakly equivalent) cylinder cofibration, 
and then take the usual colimit in the underlying category $ \cC $:
\begin{equation}
\label{colimfr}
\L{\colim}\, [\,A^{(n)} \xleftarrow{} A^{(n)}\amalg A^{(n)} \xrightarrow{} A^{(n)}] 
\,\cong\, \colim \,[\,R \xleftarrow{} R \amalg R \into \cyl(R) \,] \ .
\end{equation}
In the category of DG categories with finitely many objects, there is a canonical
cylinder object $ \cyl_{\rm BL}(R) $ defined for any semi-free DG category $ R $. We call
this object the {\it Baues-Lemaire cylinder} as it was originally constructed (in the case
of chain DG algebras) in  \cite{BL77}. The differentials of $ \cyl_{\rm BL}(R) $ are defined 
by explicit formulas in terms of differentials of $ R $, while for a semi-free resolution
$ R \sonto A^{(n)} $, the differentials of $ R $ are determined explicitly by the 
relations of $ A^{(n)} $. Thus, taking the colimit \eqref{colimfr} with  the help of
the Baues-Lemaire cylinder $ \cyl_{\rm BL}(R) $, we find an explicit presentation 
for the knot DG category $ {\mathcal A} $, given
in Definition~\ref{knot_DG_category_definition}. An elementary calculation then shows
that the DG algebra $ {\mathcal A}(1,1) $ consisting of all endomorphisms of the object `$1$'
in the DG $k$-category $ {\mathcal A} $ is precisely the knot DGA defined in \cite{Ng08}.
This explains the origin of `mysterious' algebraic formulas defining the differentials of Ng's 
combinatorial knot DGA: the differentials of the knot DGA arise from the differentials of the 
Baues-Lemaire cylinder over the canonical DG resolution of $ A^{(n)} $.



In \cite{Ng05b, Ng08}, Ng has also given an explicit description of the $0$th homology of his 
knot DGA in terms of the knot group $ \pi_1(\bR^3\!\setminus\! K) $ and the peripheral pair 
$ (m,l)$ of a meridian and longitude in $ \pi_1(\bR^3\!\setminus\! K) $.
We extend this description to the $0$th homology of the knot DG category, 
both in the framed and fully noncommutative cases 
(see Theorem~\ref{cord_category_main_theorem} and 
Theorem~\ref{noncommutative_cord_category_peripheral_description}).
Our proof is purely algebraic, in contrast to a topological proof 
given in \cite{Ng05b, Ng08}.

Finally, we mention one interesting application of our results that brings 
us back to topology. Given a link $ L \subset \bR^3 $, we consider the 
category $ \Perv(\bR^3,L) $ of perverse sheaves on $\bR^3 $ constructible with respect to the 
stratification $ L \into \bR^3 \hookleftarrow \bR^3\!\setminus\! L $ with perversity  
given by $ p(1)=0$ and $ p(3) = -1 $. Our Theorem~\ref{relation_with_perverse_sheaves}
states that $ \Perv(\bR^3,L) $ is equivalent to the category of finite-dimensional left
modules over the fully noncommutative link $k$-category  $ \tilde{A}_L$.
This leads to an algebraic description of the category $ \Perv(\bR^3,L) $ in terms of groups and quivers,
similar, in spirit, to the Gelfand-MacPherson-Vilonen description of the category 
$ \Perv(D, \{p_1, \ldots, p_n\}) $.

The paper is organized as follows.
In Section \ref{Yang_Baxter_operators_section}, we define cocartesian Yang-Baxter operators 
and the associated categorical braid closure, and give two criteria 
-- the Wada condition (Definition \ref{Wada_condition_definition}) and the Reidemeister condition (Definition \ref{Reidemeister_operators_definition}) --
for the categorical braid closure to be a link invariant (Theorem \ref{Wada_theorem_trivial_torsion} 
and Theorem \ref{Wada_theorem_nontrivial_torsion}). 
In Section \ref{Section_homotopy_braid_closure}, we extend the construction of a categorical braid closure 
to the homotopical setting. The main result in this section is Theorem \ref{invHBC}.
In Section \ref{The_GMV_operator_section}, we introduce our main example of the cocartesian 
Yang-Baxter operator associated to the GMV braid action.
In Section \ref{Section_knot_DG_cat}, we calculate the homotopy braid closure with respect to
the GMV operator, and show that the resulting DG category is an extension of Ng's knot DGA
(see Theorem \ref{explicit_DG_cat_for_homotopy_closure} and Theorem \ref{endomorphism_DGA_at_1_recovers_knot_DGA_theorem}).
The main tool in this calculation is the Baues-Lemaire cylinder on a semi-free DG category; for reader's convenience,
we review its construction is some detail.
In Section \ref{Section_cord_cat}, we compute the $0$th homology of the knot DG category,
called the {\it knot $k$-category}, and give a description  
of this category in terms of the knot group together with a peripheral pair (see Theorem \ref{cord_category_main_theorem}).
In Section \ref{FN}, we define the fully noncommutative link DG category and
extend the main results of Sections \ref{Section_knot_DG_cat} and \ref{Section_cord_cat}
to this case (see Theorem \ref{fully_noncommutative_link_DG_category_main_theorem} and Theorem \ref{noncommutative_cord_category_peripheral_description}).
While the input for the knot DG category introduced in Sections \ref{Section_knot_DG_cat} and \ref{Section_cord_cat} 
is the modified GMV action, the input for the fully noncommutative case is the original GMV action.
This allows us to relate the corresponding module category to perverse sheaves (see 
Theorem \ref{relation_with_perverse_sheaves}).
Finally, in Section \ref{Generalizations_section}, we give two natural generalizations of the 
GMV operator, inspired by the work of Wada \cite{Wad92} and Crisp-Paris \cite{CP05} in the group case.
These generalizations satisfy the Reidemeister conditions, and hence the corresponding 
homotopy braid closures give link invariants generalizing the link DG category associated to the 
original GMV action. We will describe in detail these examples in \cite{BEY}.

\vspace{1ex}

\noindent
{\bf Acknowledgements.} We would like to thank V.~Shende for bringing to our attention
and explaining to us some of his recent work on contact homology and constructible sheaves
(see \cite{ENS, NRSSZ, STZ, She, ST}). We also thank L.~Ng and B.~Cooper for interesting
questions and comments.

\section{Yang-Baxter operators and categorical braid closure}  
\label{Yang_Baxter_operators_section}

Let $\cC$ be a category closed under finite colimits. 
Let $A \in \cC$ be an object of $\cC$.  For an integer  $n \geq 2$, 
we denote the $n$-fold coproduct of copies of $A$ in $\cC$ by 
$\, A^{(n)} := A\, \amalg \stackrel{n}{\ldots} \amalg\, A\,$.
If $f : A \ra B$ is a morphism in $\cC$, we denote its $n$-fold coproduct  
by $f^{(n)} : A^{(n)} \ra B^{(n)}$.

Now, suppose that we are given an object $ A $ and a morphism $\YB $ in $\cC$.
Then, for each $\,n \geq 2$ and $ i = 1,2 \ldots, n-1$,  $\, \sigma $ induces a morphism 
$\sigma_{i,i+1}:\, A^{(n)} \to A^{(n)}$ defined by
\[ \sigma_{i,i+1} 
:= \id^{(i-1)} \, \amalg \, \sigma \, \amalg \, \id^{(n-i-1)}:\  A^{(n)} \ra A^{(n)}\]

\bdf  \label{cocartesian_Yang_Baxter_operators_definition}
A \emph{cocartesian Yang-Baxter operator} on $A$ is an invertible morphism 
$$
\YB 
$$
satisfying the equation
\begin{equation}  \label{Yang_Baxter_equation}
\sigma_{23}\,\sigma_{12}\,\sigma_{23}\, = \,\sigma_{12}\,\sigma_{23}\,\sigma_{12}  
\qquad \mbox{in}\quad \Hom_{\cC}(A^{(3)},\, A^{(3)})\ .
\end{equation}
We will often use the term ``Yang-Baxter'' as an adjective for an invertible morphism 
$\sigma$ satisfying \eqref{Yang_Baxter_equation}.
\edf

Any cocartesian Yang-Baxter operator $ \sigma $ on $A$ extends in a natural way to 
a \emph{left} action of the Artin braid group $B_n$ on $A^{(n)}$ for each $ n\geq 2$. 
We refer to this action as the action generated by $ \sigma $.

We give two basic examples of cocartesian Yang-Baxter operators corresponding to two
classical representations of the braid group $B_n$.

\begin{exmp} \label{Artin}
Let $\cC = \Gr$ be the category of groups, and let $A = \bF_1 \in \cC$ be the free group on one generator.
Consider the map $\YB$ given by
\[ 
\sigma :\, \bF_2 \ra \bF_2 \qquad x_1 \mapsto x_1 x_2 x_1^{-1}, \qquad x_2 \mapsto x_1
\]
As mentioned in the Introduction, this is a Yang-Baxter map generating the Artin 
representations\footnote{In the literature (see, e.g., \cite{Bir74}),
it is more common to extend $\sigma$ to a \emph{right} braid action.
Thus, if $\Phi : B_n \ra B_n$ is the anti-isomorphism of $B_n$ where $\Phi(\sigma_i) = \sigma_i$,
then the automorphism in the convention of \cite{Bir74} corresponding to the element $\beta \in B_n$
is equal to the automorphism in our present convention corresponding to the element $\Phi(\beta) \in B_n$.}.
\end{exmp}

\vspace{1ex}

\begin{exmp}\label{Burau}
Let $\cC = \Mod(R)$  be the category of 
modules over the commutative ring $R = \bZ[t,t^{-1}]$. Take
$A = R$ to be the free $R$-module of rank one, and define the map
$\sigma : R^{\oplus 2} \ra R^{\oplus 2}$ by left multiplication by the matrix 
$\begin{bmatrix}
1-t & t \\
1 & 0
\end{bmatrix}
$.
This map is a cocartesian Yang-Baxter operator in the category $ \Mod(R) $ 
generating the classical (unreduced) Burau representations.
\end{exmp}

\vspace{1ex}

Now, given a cocartesian Yang-Baxter operator $\YB$, we denote the resulting 
braid group action on the $n$-fold coproduct by
\[ 
\phi_n^{(A,\sigma)} : B_n \ra \Aut\,A^{(n)}
\]
Abusing notation, for a braid $\beta \in B_n$, we will often write the automorphism
$\phi_n^{(A,\sigma)}(\beta)$ simply as $\beta$ if the underlying 
Yang-Baxter operator is understood to be $(A,\sigma)$.

\bdf  \label{categorical_braid_closure_definition}
The \emph{categorical braid closure} of a braid $\beta \in B_n$ 
with respect to a cocartesian Yang-Baxter operator $\YB$ 
is defined to be the coequalizer
\[ 
\cL(A,\sigma)[\beta]\,  := \,\xcatcl{A^{(n)}}{\beta}\ ,
\]
or equivalently, the following pushout in $\cC$: 
\[
{\mathcal L}(A, \sigma)[\beta] 
= \,{{\rm colim}}\,[\,A^{(n)} \xleftarrow{\,(\beta,\, \id)\,} A^{(n)} \amalg A^{(n)} \xrightarrow{\,(\id,\, \id)\,} A^{(n)} \,]\ .
\]
\edf

\vspace{1ex}

Recall that the coequalizer of two morphisms 
$f: X \ra Y$ and $g : X \ra Y$ in $\cC$ is an object $ E \in \cC$ given together with a moprhism $p : Y \ra E$ such that the pair $(E,p)$ is universal among  all pairs such that $pf = pg$. 
In practice, computing the coequalizer amounts to taking a quotient of the object $Y$ 
by the relations $f(x) = g(x)$ for all $x \in X$.
Thus, in Example~\ref{Artin}, the categorical closure of $\beta \in B_n$ is the group presented by
\[  
\cL(A,\sigma)[\beta] = \langle \, x_1,\,\ldots\, ,\,x_n  \ | \ \beta(x_1) = x_1,\, \ldots , \,\beta(x_n) = x_n \, \rangle\ . 
\]
The next theorem is a classical result first stated by E. Artin in \cite{Art25} and proved by J. Birman in \cite{Bir74}.

\bthm[Artin-Birman]
The categorical closure of a braid $\beta \in B_n$ with respect to the Artin representation 
is the fundamental group of the link complement $\bR^3\!\setminus\! L$, 
where $L = \hat{\beta}$ is the closure of the braid $\beta$.
\ethm

Similarly, in Example~\ref{Burau}, the categorical closure of $\beta \in B_n$ 
is the module over $R = \bZ[t,t^{-1}]$
given by
\begin{equation}
\label{brb}
\cL(R,\sigma)[\beta] = 
\coker[R^{\oplus n} \xra{\id - \beta}   R^{\oplus n}] \ .
\end{equation}
In this case, we have the following theorem due to D. Goldschmidt \cite{Gol90}.

\bthm
\label{Alexm}
The categorical closure of a braid with respect to the Burau action 
is the Alexander module of the unlinked disjoint union $L \cup O$
of the braid closure $L = \hat{\beta}$ with the unknot $O$.
\ethm

Thus, the categorical braid closure of both the Artin and Burau examples 
are link invariants. This raises the natural question: When does the categorical braid closure 
of a cocartesian Yang-Baxter operator produce a link invariant?
To address this question we begin with the following definition.

Given a map $\YB$, we consider the coequalizer
\[
E := \xcoeq{A}{\sigma \circ i_2\quad}{i_2\quad}{A \amalg A}
\]
where $ i_2: A \to A \amalg A $ is the canonical map identifying
$ A $ with its second copy in $ A \amalg A $. We let 
$ p: A\amalg A \ra E$ denote the universal map such that 
$ p \circ \sigma \circ i_2 = p\circ i_2 $.

\bdf   \label{Wada_condition_definition}
We say that the map $ \sigma $ is \emph{Wada} if the following composition is an isomorphism in $ \cC $:
\[ j':\,A \xra{i_1} A\amalg A \xra{p} E \]

For a Wada map $ \sigma $, we consider the map
\[ 
j:\, A \xra{i_1} A\amalg A \xra{\sigma} A\amalg A \xra{p} E 
\]
and define the \emph{torsion} of $\sigma$ to be the map
\[ 
\chi(\sigma) = (j')^{-1} \circ j : \, A \ra A
\]
We say that a Wada map $\sigma$ has \emph{trivial torsion} if $\chi(\sigma) = \id_A $ is the identity map. 
\edf

As an easy exercise for the reader, we recommended to check that 
both the Artin and Burau Yang-Baxter operators have trivial torsion.
The next theorem explains why the categorical braid closures of these operators give link invariants.
\bthm [Wada]  \label{Wada_theorem_trivial_torsion}
Suppose that a cocartesian Yang-Baxter operator $\YB$ is Wada with trivial torsion,
then the isomorphism type of the categorical braid closure is invariant under Markov moves,
and hence give a link invariant.
\ethm

This theorem was proved in \cite{Wad92} in the special case 
when $\cC$ is the category $\Gr$ of groups, and the object $A \in \Gr$ is the 
free group $ {\mathbb F}_1 $ on one generator. However, the arguments of \cite{Wad92} can be easily
formalized and extended to a proof in the general case.

The Wada condition involves coequalizers, which makes its verification somewhat
clumsy in practice (especially, in the homotopical setting which we will discuss in the next section).
We therefore introduce another condition on $\sigma $ that, among other things, turns 
the Wada condition into a simpler form.

\bdf  \label{dualizable_definition}
We say that a map $\YB$ is \emph{dualizable} if 
\begin{enumerate}
\item The map $\YB$ is invertible.
\item The map $\sigma^R_U = (\sigma\circ i_2, i_2) :\, A\amalg A \ra A\amalg A$ is invertible.
\item The map $\sigma^L_U = (\sigma\circ i_1, i_1) :\, A\amalg A \ra A\amalg A$ is invertible.
\end{enumerate}
\edf

With this definition, we have

\bpp  
\label{torsion_for_dualizable_operators}
Assume that $\YB$ is dualizable. Consider the composition of maps
\begin{equation}  \label{j_j_maps_for_dualizable_operators}
(j,j'): \, A \amalg A \xra{ \sigma^L_U } A \amalg A  \xra{ (\sigma^R_U)^{-1} }  A \amalg A  \xra{\nabla} A \ ,
\end{equation} 
where $ \nabla $ is the canonical folding map.
Then, $\sigma$ is Wada if and only if the map $j'$ is an isomorphism.
In this case, the torsion of $\sigma$ is given by
\begin{equation}  \label{dualizable_torsion_formula}
\chi(\sigma) = (j')^{-1} \circ j
\end{equation} 
\epp

\vspace{1ex}

Next, we introduce our main definition.
\bdf  
\label{Reidemeister_operators_definition}
We say a map $\YB$ is \emph{Reidemeister} if it is Yang-Baxter, dualizable, 
and Wada with invertible torsion $\chi(\sigma) $ (see 
Proposition \ref{torsion_for_dualizable_operators}). 
\edf

In Theorem \ref{Wada_theorem_trivial_torsion}, we required the torsion of $\sigma$ to be trivial.
However, in our main example of the cocartesian Yang-Baxter operator associated to the GMV action
(see Section \ref{The_GMV_operator_section},) the torsion is not trivial, but invertible.
It turns out, if $\sigma$ is Reidemiester, with not necessarily trivial torsion, then an analogue 
of Theorem \ref{Wada_theorem_trivial_torsion} holds, provided we modify the 
categorical braid closure in an appropriate way.
From now on, for simplicity, we will work with knots ({\it i.e.}, links with one component), while all that follows 
holds for the general case of links (See \cite{BEY} for more details).

Let $\YB$ be a Reidemeister operator with (invertible) torsion $\chi : A \ra A$.
Suppose that $\beta \in B_n$ is a braid that closes to a knot $\hat{\beta} = K$.
\bdf  \label{writhe_adjusted_categorical_braid_closure_knots_definition}
The \emph{writhe-adjusted categorical closure} of  $\beta$ 
with respect to $\sigma$ is defined by
\[ \cL^{\wa}(A,\sigma)[\beta] = \xcatcl{A^{(n)}}{ \Psi_0 } \]
The map $\Psi_0 :\, A^{(n)} \ra A^{(n)} $ is the composition of morphisms
\[ \Psi_0 : A^{(n)} \xra{ \chi^{-w} \, \amalg \, \id^{(n-1)} }   A^{(n)}  \xra{\beta}   A^{(n)}\]
where $w=w(\beta)$ is the writhe of the braid $\beta$.
\edf

We have
\bthm  
\label{Wada_theorem_nontrivial_torsion}
For any Reidemeister operator $\YB$, 
the isomorphism type of the writhe-adjusted categorical braid closure
is invariant under Markov moves, and hence gives a knot invariant.
\ethm

There is a conceptual way to prove Theorems \ref{Wada_theorem_trivial_torsion} 
and \ref{Wada_theorem_nontrivial_torsion} by interpreting the categorical
braid closure as an abstract trace in the sense of \cite{JSV96}.
To this end, we extend the category $\cC$ to a larger category $\hatC$ with the same object set 
$\Ob(\hatC) = \Ob(\cC)$.
The hom-set $\Hom_{\hatC}(A,B)$ in $\hatC$ is given by the set of cospans 
$[B \ra X \la A]$ modulo isomorphisms that are identity on $A$ and $B$.
The composition in $\hatC$ is defined by pushouts in the obvious way.
The category $\hatC$ has a monoidal product $\boxtimes$ induced by the coproduct in $\cC$;
hence the monoidal structure of $ \hatC $ extends the cocartesian monoidal structure on $\cC$
in the sense that there is a faithful, strongly monoidal functor $\iota : (\cC, \amalg) \into (\hatC , \boxtimes)$.
Moreover, the monoidal category $(\hatC , \boxtimes)$ has a canonical pivotal structure,
and therefore an abstract trace axiomatized in \cite{JSV96} and \cite{Sel}. 
%
%
%
This construction allows us to interpret the categorical braid closure as an abstract trace in $\hatC$, 
and the Wada condition (Definition \ref{Wada_condition_definition}) can then be interpreted as a 
condition on the partial trace of $\sigma$ viewed as a morphism in $\hatC$ under the faithful embedding $\iota : \cC \into 
\hatC$. We remark that the Wada condition reinterpreted this way is analogous to a condition on partial trace 
for ``enriched Yang-Baxter operators'' introduced by Turaev in \cite{Tur88}.
This interpretation  allows us to prove Theorem \ref{Wada_theorem_trivial_torsion} 
and Theorem \ref{Wada_theorem_nontrivial_torsion} by diagrammatic tensor calculus.
In fact, starting with a Reidemeister operator, one can construct a ribbon category 
(in the sense of \cite{RT90, Tur10}, see also \cite{JS91a}, \cite{Shu94}), whose associated 
link invariant, which lives in the set $\Hom_{\hatC}(\phi, \phi)$ of isomorphism classes 
of objects in $\cC$, coincides with the categorical braid closure.
One can see this as another justification for the term ``categorical braid closure''.
For details, we refer the reader to \cite{BEY}.

\brm   
\label{relation_with_quandles_remark}
The notion of a cocartesian Yang-Baxter operator is closely related to biracks and biquandles\footnote{For 
the definition and basic examples of biracks and biquandles we refer to \cite{FJK04, CSWES09}.}.
To make this relation precise we recall the (dual) Yoneda embedding $ \cC \into {\rm Fun}(\cC, \Set) $
for a category $ \cC $ that associates to an object $A \in \cC $ the corepresentable functor 
$\,h^A := \Hom(A,\,\mbox{--}\,):\, \cC \ra \Set\,$. A cobirack structure on $A$ can then be defined by factoring 
$h^A$ into a composition of functors $\cC \ra \Birack \xra{\text{forget}} \Set$. 
Similarly, a cobiquandle structure on $A$ is the factorization of $h^A$ into a composition of 
functors $\cC \ra \Biquandle \xra{\text{forget}} \Set$.
Then, one can show that, giving a cobirack structure on $A$ is equivalent to giving 
a dualizable cocarteisan Yang-Baxter operator on $A$.
Similarly, giving a cobiquandle structure on $A$ is equivalent to giving 
a Reidemeister operator on $A$ with trivial torsion.
Thus, in particular, given a Reidemeister operator on $A$ with trivial torsion, the set $X_B = \Hom_{\cC}(A,B)$ 
has a natural biquandle structure for any object $B \in \cC$.
In fact, many examples of biquandles in the literature arise in this manner. 
(In particular, almost all examples of biquandles given in \cite{FJK04} are of this form.)

Given a biquandle $X$ and a link $L$, one can define a combinatorial link invariant 
called $\Color_X(L)$, which is the set of colorings of a link diagram of $L$ by the biquandle $X$
(see \cite{FJK04}). If the link $L$ is the closure of a braid $\beta \in B_n$,
and if the biquandle $X = X_B$ arises from a cobiquandle structure on an object $A$
in the sense above, then we have 
\[ 
\Color_{X_B}(L)   = \Hom_{\cC}( \cL(A,\sigma)[\beta] , B ) \ .  
\]
This last formula gives a combinatorial interpretation of the categorical braid closure in terms of 
arcs of a link diagram that we alluded to in the Introduction.

\erm

\section{Homotopy braid closure}  
\label{Section_homotopy_braid_closure}
In this section, we will work with model categories and assume the reader to have 
some familiarity with the theory of model categories and derived functors.
For an excellent introduction, we recommend the Dwyer-Spalinski article \cite{DS95}, 
which covers enough material  for the present paper. 
For a more comprehensive study of model categories, we refer to \cite{Hir03, Hov99}.

If $\cC$ is a model category,
then the notion of a categorical braid closure associated to a cocartesian Yang-Baxter operator in 
$\cC $ admits a natural generalization, which is obtained by replacing colimits in 
Definition \ref{categorical_braid_closure_definition} 
and Definition \ref{writhe_adjusted_categorical_braid_closure_knots_definition}
by homotopy colimits.
The analogues of Theorem \ref{Wada_theorem_trivial_torsion} and 
Theorem \ref{Wada_theorem_nontrivial_torsion}  hold in this homotopical setting,
provided that the object $A$ satisfies a pseudoflatness condition (which roughly says 
that the $n$-fold homotopy coproduct of $A$ coincides with the $n$-fold coproduct of $A$).
Since the translation to the homotopical context is faily straightforward,
we omit here formal details. Instead, we will give explicit definitions and statements only in
the case, where the notion of a writhe-adjusted braid closure
is further refined by allowing the knot in question to be colored by certain maps.

\bdf  \label{sigma_natural_definition}
Given a cocartesian Yang-Baxter operator $\YB$, we
say that a map $\theta : A \ra A$ is $\sigma$-\emph{natural} if 
the following two diagrams commute
\[ \xymatrix{
A \amalg A \ar[r]^{\sigma} \ar[d]_{\theta \, \amalg \, \id}  & A \amalg A  \ar[d]^{\id \, \amalg \, \theta} 
& & A \amalg A \ar[r]^{\sigma} \ar[d]_{\id \, \amalg \, \theta}  & A \amalg A  \ar[d]^{\theta \, \amalg \, \id} \\
A \amalg A \ar[r]^{\sigma}  & A \amalg A   & & A \sqcup A \ar[r]^{\sigma}  & A \amalg A
}\]
\edf

Now, let $\YB$ be a Reidemeister operator with (invertible) torsion $\chi : A \ra A$,
and let $\theta : A \ra A$ be a $\sigma$-natural map.
Suppose that $\beta \in B_n$ is a braid that closes to a knot $\hat{\beta} = K$.
\bdf  \label{writhe_adjusted_colored_homotopy_braid_closure_knots_definition}
The \emph{colored, writhe-adjusted homotopy closure} of  $\beta$ 
with respect to $\sigma$ and the coloring $\theta$ 
is defined to be the homotopy coequalizer
\[ \bm{L} \cL^{\wa}(A,\sigma)[\beta] = \xhocl{A^{(n)}}{ \Psi }\]
where $\Psi : A^{(n)} \ra A^{(n)}$ is the composition
\begin{equation}   \label{Psi_map_definition}
A^{(n)} \xra{ \theta \chi^{-w} \, \amalg \, \id^{(n-1)} } A^{(n)} \xra{\beta} A^{(n)} 
\end{equation}
\edf
When $\theta = \id_A$ is the identity map, which is always $\sigma$-natural, this construction 
obviously reduces to the original one without coloring.

\bthm
\label{invHBC}
Let $\YB$ be a Reidemeister operator on a pseudoflat object $A$ in a model category $\cC$,
and let $\theta : A \ra A$ is a $\sigma$-natural map.
Then, the isomorphism type in $\HoC$ $($i.e., the weak equivalence type in $\cC$ $)$
of the colored, writhe-adjusted homotopy braid closure
is invariant under Markov moves, and hence gives a knot invariant.
\ethm

This theorem allows one to refine many classical link invariants defined by categorical 
braid closure. To illustrate this we will return to Examples \ref{Artin} and \ref{Burau} 
in Section~\ref{Yang_Baxter_operators_section}.

\begin{exmp}
To compute the homotopy braid closure of the Burau representations (see Example \ref{Burau}),
we embed the module category $ \Mod(R) $ into the category  $ \Ch(R)$ of chain complexes  
in the usual way. The module $A = R$ is then identified with the chain complex 
$ A = [0 \to R \to 0] $, with $R$ concentrated in degree $0$. 
The category $ \Ch(R) $ has a natural (projective) model structure, 
with weak equivalences being the quasi-isomorphisms and the fibrations
being the degreewise surjective morphisms of complexes (see \cite[Section~2.3]{Hov99}). 
Every object in this model category is pseudoflat. The corresponding homotopy category 
$ \Ho(\cC) $ is the (unbounded) derived category $ {\mathcal D}(R) $ of $R$-modules. Now, 
the homotopy closure of a braid $\beta \in B_n$ with respect to the Burau operator 
$ \sigma: R^{\oplus 2} \to R^{\oplus 2} $ is given by the mapping cone of the morphism
$\, \id - \beta \,$ in the derived category $ {\mathcal D}(R) $: {\it i.e.},
$$
\L\cL(R,\sigma)[\beta] = {\rm C\hat{o}ne}(\id - \beta)\, := \, 
[\,0 \ra R^{\oplus n} \xra{\id - \beta} R^{\oplus n} \ra 0\,]\ ,
$$
where the two copies of $ R^{\oplus n} $ are concentrated in homological degrees $0$ and $1$.
The isomorphism class of $ \L\cL(R,\sigma)[\beta] $ in $ {\mathcal D}(R) $ is a link invariant 
by Theorem \ref{invHBC}. Note that the homology of this complex in degree $0$ 
is the cokernel of $\, \id - \beta \,$, which is precisely
the categorical braid closure $ \cL(R,\sigma)[\beta] $, see \eqref{brb}.
\end{exmp}

\vspace{1ex}

\begin{exmp}
To compute the homotopy braid closure of the Artin representations (see Example \ref{Artin}),
we embed the category  of groups into the model category $\sGr$ of {\it simplicial} groups.
The model structure on $\sGr$ is inherited from the category $ \mathtt{sSet} $ of simplicial sets, 
so that the weak equivalences and fibrations of simplicial groups are the weak equivalences and 
fibrations of the underlying simplicial sets  (see, e.g., \cite{GJ09}). 
The category of simplicial groups has a rich homotopy 
theory which is classically known to be equivalent to that of topological spaces.
Precisely, there is a Quillen equivalence between $\sGr $ and 
the category $ \mathtt{Top}_{0, \ast} $ of connected pointed topological spaces given 
by the composition of functors\:
\begin{equation}
\label{KMeq}
\sGr \xrightarrow{\overline{W}} \mathtt{sSet} \xrightarrow{|\,\mbox{--}\,|} \mathtt{Top}_{0,\ast} 
\end{equation}
where $ \overline{W} $ is Kan's bar construction assigning to a simplicial group its classifying
simplicial space and $\,|\,\mbox{--}\,| \,$ is Milnor's geometric realization functor.
The functor \eqref{KMeq} induces an equivalence of the homotopy categories $\,
\Ho(\sGr) \cong \Ho(\mathtt{Top}_{0,\ast})\,$, which gives a bijective correspondence
between the homotopy classes of simplicial groups and the 
homotopy classes of pointed connected CW complexes (see, e.g., \cite[V.6.4]{GJ09}).
In this way, every simplicial group can be thought of as representing  a topological space (up to homotopy).

Now, if we regard  $\bF_1$ as a discrete simplicial group in $ \sGr $, then
the homotopy closure $ \L\cL(\bF_1,\sigma)[\beta] $ of a braid $ \beta \in B_n$ with respect to 
the Artin operator $ \sigma: \bF_2 \to \bF_2 $ is represented by a simplicial group in $ \Ho(\sGr) $ that 
corresponds under the above equivalence to the link complement $\bR^3\! \setminus\! L$ of the closure 
$L = \hat{\beta}\,$: {\it i.e.},
$$
|\,\overline{W}\,\L\cL(\bF_1,\sigma)[\beta]\,| \,\simeq\, \bR^3\! \setminus\! L\ .
$$
Thus, in the case of Artin representations, the homotopy braid closure completely recovers the
homotopy type of the space $\bR^3 \setminus L$, while the categorical braid closure gives only its
fundamental group (see \cite{BEY}).
\end{exmp}

\brm
Theorem~\ref{invHBC} holds in a more general case, when 
various conditions on the map $\sigma $ hold only up to homotopy.
For example, one can require the Yang-Baxter equation \eqref{Yang_Baxter_equation} 
to hold only in $\HoC$, the maps in Definition \ref{dualizable_definition} to be only weak equivalences
and the map $j : A \ra A$ defined in \eqref{j_j_maps_for_dualizable_operators}
to be only an isomorphism in $\HoC$. The braid group action, as well as the torsion map, are then 
only defined in $\HoC$. This makes a precise definition of a homotopy braid closure a little tedious. 
We omit it here, referring the reader to \cite{BEY} instead.
\erm

\section{The Gelfand-MacPherson-Vilonen action}   
\label{The_GMV_operator_section}
In this section, we fix a commutative ring $k$ with unit.
By a {\it $k$-category}, we mean a category enriched over the category of $k$-modules.
Let $\Cat^*_k$ be the category of all (small) pointed $k$-categories,
where a $k$-category $\mathscr{A}$ is {\it pointed} if there is a distinguished object $* \in \mathscr{A}$.
Maps ({\it i.e.}, $k$-linear functors) between pointed $k$-categories are required to preserve the distinguished objects.

As in the Introduction, we consider the path category $k\langle Q \rangle$ of the quiver 
$\, Q = \bigl[  
\vcenter{
\xymatrix@1{
\infinity  \ar@/^/[rr]^{a}  
&& 0  \ar@/^/[ll]^{a^{*}} 
}
} \bigr]
\,$
Let $T \in k\langle Q \rangle(0,0)$ be the element in the endomorphism algebra 
of $k\langle Q \rangle$ of 
the object $0$ defined by $T = e_0 + a a^*$, and let $\tilde{A}$ be the $k$-category 
$ \tilde{A} = k\langle Q \rangle [T^{-1}]$,
which is pointed by taking the object $0$ as the distinguished object.

The coproduct in  $\Cat^*_k$ is given by the \emph{fusion product}.
{\it i.e.} the coproduct of $X, Y \in \Cat^*_k$ is the $k$-category obtained by collapsing the two distinguished objects in the disjoint union of $X$ and $Y$ into a single object.
In particular, the $n$-fold coproduct of $\tilde{A} \in \Cat^*_k$ is the $k$-category $\tilde{A}^{(n)}$ defined 
in \eqref{widetilde_A}.

We begin with the following result mentioned in the Introduction.

\bthm[Gelfand, MacPherson, Vilonen]  \label{GMV_is_Yang_Baxter_theorem}
The map $\sigma : \tilde{A}^{(2)} \ra \tilde{A}^{(2)}$ 
defined by \eqref{GMV1} is a cocartesian Yang-Baxter operator on the object $\tilde{A}$ 
in the category $\Cat^*_k$. We call $ \sigma $ the GMV operator.
\ethm
The GMV operator induces an action of $B_n$ on $\tilde{A}^{(n)}$,
where the generator $\sigma_i \in B_n$ acts on objects by swapping $i$ and $i+1$,
while fixing all other objects, and on morphisms by formula \eqref{GMV}.
The next observation is straightforward to check.

\begin{lem}  
\label{GMV_operator_Reidemeister}
The GMV operator \eqref{GMV1} is Reidemeister with torsion given by
\begin{equation} 
\label{GMV_operator_torsion}
\chi : \tilde{A} \ra \tilde{A}, \qquad a \mapsto Ta, \qquad a^* \mapsto a^* T^{-1}
\end{equation}
\end{lem}

As explained in the Introduction, formula \eqref{GMV} for the braid group action first 
appeared in \cite{GMV96} in relation to perverse sheaves. More precisely, it was shown 
in \cite{GMV96} that any choice of `cuts' 
({\it i.e.}, a family $\Theta$ of $n$ simple curves on $D\setminus\{p_1,\ldots,p_n\}$, 
going from a chosen point near $p_i$ to the chosen endpoint $p_0$ near the boundary $\partial D$,
so that any two such curves intersect only at $p_0$)
induces an equivalence of categories
$\,\tilde{E}_{\Theta} : \Perv(D, \{p_1, \ldots ,p_n\}) \simeq \tilde{\mathscr{Q}}\,$
from the category of perverse sheaves on the disk $D$ 
with only possible singularities at the points $\{p_1,\ldots,p_n\}$,
to a quiver category $\tilde{\mathscr{Q}}$ \emph{isomorphic}
to the category $\Mod(\tilde{A}^{(n)})$ of finite-dimensional 
modules over the $k$-category $\tilde{A}^{(n)}$.

Now, the braid group $B_n$ acts as a mapping class group on the disk $D$ with
$n$ marked points $ \{p_1, \ldots ,p_n\} $, and hence
acts (in a certain sense) on the category  $\Perv(D, \{p_1, \ldots ,p_n\})$. 
If we fix a family $\Theta $ of cuts, 
this translates to an action of $B_n$ on the quiver category $\tilde{\mathscr{Q}}$.
In fact, it is shown in \cite{GMV96} that there is a \emph{strict} action of $B_n$ on 
the quiver category $\tilde{\mathscr{Q}}$ that coincides under the equivalence $\tilde{E}_{\Theta}$ 
with the natural action on the category $\Perv(D, \{p_1, \ldots ,p_n\})$ (see \cite[Proposition 1.3]{GMV96}).
This strict $B_n$ action on the quiver category $\tilde{\mathscr{Q}}$
is in fact induced by the action \eqref{GMV} on the $k$-category $\tilde{A}^{(n)}$.
More precisely, the left action \eqref{GMV} induces a strict left action on the module category $\Mod(\tilde{A}^{(n)})$ 
where $\beta \in B_n$ acts by $M \mapsto (\beta^{-1})^*(M)$.
This coincides under the isomorphism of categories $\Mod(\tilde{A}^{(n)}) \cong \tilde{\mathscr{Q}}$
with the strict $B_n$ action on 
the quiver category $\tilde{\mathscr{Q}}$ constructed in \cite{GMV96}.

\brm
In the notation of \cite{GMV96}, a module $M \in \Mod(\tilde{A}^{(n)})$ corresponds to the representation
of the quiver $ Q^{(n)} $ which is given by a collection of vector spaces and maps $M(0) = A$, $M(i) = B_i$, $M(a_i) = q_i$, $M(a_i^*) = p_i$. The functors $ (\sigma_i^{-1})^* $ on modules correspond to the operations denoted by $T_i$. For example, $( (\sigma_i^{-1})^*M )(a_j)$  means 
$ T_i(q_j)$ in the notation of \cite{GMV96} .
\erm

Next, we introduce a slight modification of the GMV action. Let $\overline{Q}_n$ denote the following quiver 
\begin{equation}   
\label{quiver_overline_Qn_definition}
\overline{Q}_n := \qquad 
\vcenter{
\xymatrix{
1  \ar@/^3pc/[rrrr]^{a_1} \ar@/_2pc/[rrrr]^{a_n}  \ar@{}[rrrr]|{\vdots} & & & & 
0  \ar@/_2pc/[llll]^{a_1^*}   \ar@/^3pc/[llll]^{a_n^*}
}
}
\end{equation}
Fix an invertible element $\mu \in k^{\times}$, and define  the $k$-category $A^{(n)}$ by
\begin{equation} 
 \label{An_definition}
A^{(n)} = k \langle \overline{Q}_n \rangle / (a_i^* a_i = (\mu - 1) e_1)_{i=1,\ldots,n}
\end{equation}
Notice that the elements $\,T_i = e_0 + a_i a_i^* $ are invertible in $A^{(n)}$
for all $i=1,2, \ldots, n$. Hence, formula \eqref{GMV} still defines 
a braid group action on $A^{(n)}$. 

The $k$-category $A^{(n)}$ is obtained from $\tilde{A}^{(n)}$ by applying the following two operations:
\begin{enumerate}
\item taking the quotient of $ \tilde{A}^{(n)}$ modulo the relations $a_i^* a_i = (\mu - 1) e_i$,
\item collapsing the vertices $ \{1,\ldots,n\} $ into a single vertex $1$.
\end{enumerate}
The GMV braid action on $\tilde{A}^{(n)}$ descends to a braid action on $A^{(n)}$, which we will call
the \emph{$\mu$-central GMV action}.

One advantage of working with the $\mu$-central GMV action is that it
fixes the set of objects of $A^{(n)}$. In particular, one can consider the
 induced braid  action on the endomorphism algebra of any object of $A^{(n)}$. 
 Specifically, let $ A^{(n)}(1,1) $ denote the endomorphism algebra of the object `$1$'
 in $ A^{(n)}$. For $ i,j=1,2, \ldots, n$, consider the elements 
 $\, A_{ij} := - a_i^* a_j \in  A^{(n)}(1, 1) $. Then, it is easy to see  that
 the  algebra $A^{(n)}(\infinity, \infinity)$ has the the following presentation 
\[ 
A^{(n)}(\infinity, \infinity) = k \langle  A_{ij} \rangle / (A_{ii} = 1 - \mu )
\]
It is straightforward to compute the induced braid group action on this  algebra
in terms of the generators $ A_{ij} $. However, we will write the corresponding formulas
in terms of other generators $ a_{ij} $ related to  $ A_{ij} $ by a simple rescaling:
\begin{equation}  
\label{Aij_aij_relation_substitution}
a_{ij} := 
\begin{cases}
A_{ij}\ , &i<j \\
-\mu^{-1} A_{ij}\ , &i>j
\end{cases}
\end{equation}
The associative algebra $A^{(n)}(\infinity, \infinity)$ is free on these generators,
and the braid group action on $A^{(n)}(\infinity, \infinity)$ is given by
\begin{equation} 
\label{Magnus_action_equation}
\sigma_k \, :  \, \begin{cases}
a_{ki} \, \mapsto \, a_{k+1,i} - a_{k+1, k} \,  a_{ki} &\ (i \neq k, k+1)\\
a_{ik} \, \mapsto \, a_{i, k+1} - a_{ik} \, a_{k,k+1} &\ (i \neq k, k+1)\\
a_{k+1,i} \, \mapsto \, a_{ki}  &\ (i \neq k, k+1) \\
a_{i,k+1} \, \mapsto \, a_{ik} &\  (i \neq k, k+1)\\
a_{k, k+1}  \, \mapsto \, - a_{k+1, k}  &  \\
a_{k+1, k}  \, \mapsto \, - a_{k, k+1} & \\
a_{ij} \, \mapsto \, a_{ij}  &\  (i,j \neq k, k+1 )
 \end{cases}
\end{equation}
Formulas \eqref{Magnus_action_equation} first appeared in \cite{Hum92, Hum01} as a generalization of the classical Magnus action \cite{Mag80}; we therefore call \eqref{Magnus_action_equation} the \emph{Humphries-Magnus braid action}. The Humphries-Magnus braid action was used by Ng in \cite{Ng05a, Ng05b, Ng08, Ng11, Ng14} as part of his definition
of the (combinatorial) knot DGA (see, e.g., \cite{Ng14}, Definition 3.3). Now,
the main results of the present paper in relation to Ng's work can be summarized schematically by the following diagram

\begin{equation}  \label{schematic_diagram_square_root_of_DGA}
\xymatrix{ *+[F]{\txt{$\mu$-central\\GMV action}} 
\ar[rrr]^{A^{(n)} \, \mapsto \, A^{(n)}(\infinity,\infinity)}
\ar[dd]_{\txt{Homotopy\\braid closure}}
& & & *+[F]{\txt{Humphries-Magnus\\action}} \ar[dd]^{\txt{Ng's construction}}   \\
\\
*+[F]{\txt{Knot DG category}}
\ar[rrr]^{\cA \, \mapsto \, \cA(\infinity, \infinity)}
& & & *+[F]{\txt{Knot DGA}}  }
\end{equation} 

\section{The knot DG category}  
\label{Section_knot_DG_cat}
Let $\cC = \dgcat^{ \{ 0 , 1 \} }_k$ be the category comprising all small DG $k$-categories 
with object set $\{0,1\}$. The morphisms of such DG categories in $\cC$ are required to be 
the identity map on the object set $\{0,1\}$. The $k$-category $A^{(n)}$ defined in \eqref{An_definition} can then be identified with the $n$-fold coproduct of copies of $A := A^{(1)}$ in the category $\cC$.
Moreover, the $\mu$-central GMV action is induced by a cocartesian Yang-Baxter operator $\YB$ 
given by the same formula as in \eqref{GMV1}.
The same calculation as in Lemma \ref{GMV_operator_Reidemeister} shows 
that this cocartesian Yang-Baxter operator is Reidemeister with torsion given by  formula \eqref{GMV_operator_torsion}.
Moreover, the following lemma gives a $\sigma$-natural map (Definition \ref{sigma_natural_definition}) that can be used to color a knot.
\begin{lem}
For any invertible element $\lambda \in k^{\times}$, the map
$ \theta_{\lambda} : A \ra A$ given by $ (a,a^*) \mapsto ( \lambda^{-1} a, \lambda a^*) $
is $\sigma$-natural.
\end{lem}
The category $\cC = \dgcat^{ \{ 0 , 1 \} }_k$ has a model structure, in which 
a morphism $f : X \ra Y$ is a weak equivalence (resp., fibration) if and only if for any pair
of objects $a,b \in X$, the map $f : X(a,b) \ra Y(a,b)$ is a quasi-isomorphism (resp., surjection) 
of chain complexes. This model category is cofibrantly generated, therefore the
cofibrations can be characterized as retracts of relative cell complexes (see \cite{Hir03}, \cite{BEY}),
which in particular, include semi-free extensions by arrows in non-negative (homological) degree.

One can show that the $k$-category $A$ viewed as an object of the model category
$\cC = \dgcat^{ \{ 0 , 1 \} }_k$ is pseudoflat. Therefore, the colored writhe-adjusted homotopy braid closure (Definition \ref{writhe_adjusted_colored_homotopy_braid_closure_knots_definition}) with respect to the GMV operator $(A,\sigma)$ and the coloring $\theta_{\lambda}$:
\begin{equation} \label{homotopy_closure_of_GMV_hocoeq_form}
\begin{split}
\bm{L}\cL^{\wa}(A,\sigma)[\beta; \theta_{\lambda}] 
& := \xhocl{A^{(n)}}{\Psi} \\
& = \bm{L} \colim \bigl[ \, A^{(n)} \xla{(\Psi, \id)} A^{(n)} \amalg A^{(n)} \xra{(\id,\id)} A^{(n)} \, \bigr]
\end{split}
\end{equation}
gives a quasi-isomorphism type in the category $\cC = \dgcat^{ \{ 0 , 1 \} }_k$,
which is a knot invariant. 

We now describe this knot invariant in explicit terms.
Let $\cQ$ be the following \emph{graded} quiver
\[
\cQ = \qquad 
\vcenter{
\xymatrix{
\infinity  
\ar@/^4.5pc/[rrrr]|{b_1}   \ar @{{}{ }{}}@/^4pc/[rrrr]|{\vdots}   \ar@/^3pc/[rrrr]|{b_n} 
\ar@/^2pc/[rrrr]|{a_1}   \ar @{{}{ }{}}@/^1.5pc/[rrrr]|{\vdots}   \ar@/^0.5pc/[rrrr]|{a_n}  
\ar@`{(-20,20),(-20,-20)}[]|{\eta_n}
\ar@{{}{ }{}}@`{(-25,25),(-25,-25)}[]|{\cdots}
\ar@`{(-30,30),(-30,-30)}[]|{\eta_1}
 & & & & 
0  
\ar@/^4.5pc/[llll]|{b_1^*}   \ar @{{}{ }{}}@/^3.5pc/[llll]|{\vdots}   \ar@/^3pc/[llll]|{b_n^*} 
\ar@/^2pc/[llll]|{a_1^*}   \ar @{{}{ }{}}@/^1pc/[llll]|{\vdots}   \ar@/^0.5pc/[llll]|{a_n^*}  
}
}
\]
where the degrees of arrows are assigned by
\[
\begin{split}
\deg(a_1) &= \ldots = \deg(a_n) = \deg(a_1^*) = \ldots = \deg(a_n^*) = 0 \\
\deg(b_1) &= \ldots = \deg(b_n) = \deg(b_1^*) = \ldots = \deg(b_n^*) = 1 \\
\deg(\eta_1) &= \ldots = \deg(\eta_n) = 2
\end{split}
\]
Let $\beta \in B_n$ be a braid that closes to a knot $K$. 

\bdf   
\label{knot_DG_category_definition}
We defined the \emph{knot DG category of $K$} to be the DG $k$-category
\[ 
\cA_K = k\langle \cQ \rangle/(a_i^* a_i = (\mu-1)e_{\infinity})_{1\leq i \leq n} 
 \]
with differential given by
\begin{equation}\label{difform}
\begin{split}
 d(b_i) &=  \Psi(a_i) - a_i  \\
 d(b^*_i) &= \Psi(a^*_i) - a^*_i \\
 d(\eta_i) &= - b^*_i a_i - \Psi(a^*_i) b_i 
\end{split}
\end{equation}
where $\Psi : A^{(n)} \ra A^{(n)}$ is the map defined in \eqref{Psi_map_definition}.
\edf

Our main results regarding the knot DG category $ \cA = \cA_K $ can be encapsulated into the following two theorems.
\bthm   \label{explicit_DG_cat_for_homotopy_closure}
The knot DG category $\cA$ represents the quasi-isomorphism type of the 
homotopy coequalizer \eqref{homotopy_closure_of_GMV_hocoeq_form}.
Therefore, the quasi-isomorphism type of the knot DG category $\cA$
is independent of the choice of a braid $\beta$ that closes to a given knot $K$. 
\ethm
\bthm  \label{endomorphism_DGA_at_1_recovers_knot_DGA_theorem}
Let the base commutative ring be $k = \mathbb{Z}[\mu^{\pm 1}, \lambda^{\pm 1}]$ Then
the quasi-isomorphism type of the endomorphism DG algebra $\cA(1,1)$
coincides with the quasi-isomorphism type of the knot DGA constructed in \cite{Ng08}.
\ethm

Theorem~\ref{endomorphism_DGA_at_1_recovers_knot_DGA_theorem} gives an alternative proof
of one of the main results in \cite{Ng05a, Ng08} that states that the underlying quasi-isomorphism type of the combinatorial knot DGA  is a knot invariant\footnote{Note, however, that the results in \textit{loc. cit.} are slightly stronger as they refer to the invariance of the stable tame isomorphism type rather than the quasi-isomorphism type of the corresponding knot DGA.}.

\begin{proof}[Proof of Theorem \ref{endomorphism_DGA_at_1_recovers_knot_DGA_theorem}]
Define the following morphisms in $\cA$, which are elements of the endmorphism DG algebra
$ \cA(1,1) $ of different homological degrees:
\begin{equation}  \label{generators_ABCDe_equation}
\begin{split}
 A_{ij} &= -a_i^{*} a_j \in \cA(\infinity,\infinity)_0 \\
 B_{ij} &= b_i^{*} a_j \in \cA(\infinity,\infinity)_1 \\
 C_{ij} &= a_i^{*} b_j \in \cA(\infinity,\infinity)_1 \\
 D_{ij} &= b_i^{*} b_j \in \cA(\infinity,\infinity)_2 \\
 e_i &= -\eta_i  \in \cA(\infinity,\infinity)_2 
\end{split}
\end{equation}
Then, the DG algebra $\cA(\infinity,\infinity)$ is freely generated by the elements 
\eqref{generators_ABCDe_equation}, modulo the relations $A_{ii} = 1-\mu$.
The differentials of these elements can be easily computed by the Leibnitz rule,
using formulas \eqref{difform}:
\begin{equation} \label{differentials_for_endomorphism_DGA_first_form}
\begin{split}
d(A_{ij}) &= 0  \\
d(B_{ij}) &= \Psi(a_i^*)a_j - a_i^* a_j \\
d(C_{ij}) &= a_i^* \Psi(a_j) - a_i^* a_j \\
d(D_{ij}) &= (\Psi(a_i^*) - a_i^*)b_j + b_i^*(\Psi(a_j) - a_j) \\
d(e_i) &= b_i^* a_i + \Psi(a_i^*)b_i
\end{split}
\end{equation}
This explicit description allows one to identify $\cA(\infinity,\infinity)$ with
the combinatorial knot DGA as defined in \cite[Definition~2.6]{Ng08} (see also \cite{Ng14}). 
See \cite{BEY} for details of this calculation.
\end{proof}

To prove Theorem \ref{explicit_DG_cat_for_homotopy_closure} one has to calculate the 
homotopy pushout \eqref{homotopy_closure_of_GMV_hocoeq_form}.
As shown in \cite{BEY}, it suffices for this to resolve 
the right-pointing arrow by a strong cofibration 
({\it i.e.}, a cofibration whose domain is cofibrant), 
and then take the ordinary pushout of the resulting diagram.
Thus, we need to find a semi-free resolution $p : B \sonto A$ and then construct
an appropriate cylinder object $\cyl(B)$ on $B$.
The right-pointing arrow in the pushout diagram in
\eqref{homotopy_closure_of_GMV_hocoeq_form}
will then be resolved by taking the $n$-fold coproduct $\cyl(B)^{(n)}$ of this cylinder object.

To construct a semi-free resolution of $A$, we consider
the \emph{graded} quiver
\begin{equation}  \label{quiver_widetilde_Q_definition}
\tilde{Q} = \qquad 
\left[  \qquad
\vcenter{
\xymatrix{
\infinity \ar@`{(-10,10),(-10,-10)}[]_{\xi} \ar@/^/[rr]^{a}  
&& 0  \ar@/^/[ll]^{a^{*}} 
}
} 
\qquad \right]
\end{equation} 
where $\deg(a) = \deg(a^*)=0$ and $\deg(\xi) = 1$.
Define $B \in \cC$ to be the semi-free DG category
$\, B := k \langle \tilde{Q} \rangle\,$ with differential 
given by $d \xi = (a^{*}a - (\mu-1)e_{\infinity})$.
Then, one can show (see \cite{BEY}) that the canonical map
\[ p : B \ra A, \qquad a \mapsto a, \quad a^* \mapsto a^*, \qquad \xi \mapsto 0 \] 
is a quasi-isomorphism. Thus, $ B $ can be used as a cofibrant replacement for $A$.

Next, to define a cylinder on $B$ we will use the construction of canonical cylinder objects 
for semi-free DG algebras given in \cite{BL77}. This construction will play an important role
in our calculations, so we review it in some detail.

Let $R$ be a DG algebra whose underlying graded algebra is free over a graded $k$-module $V$.
We write $R = (T(V),d)$.
Let $ \cyl(R) $ be the graded algebra defined by $ \cyl(R) := T(V \oplus V' \oplus sV)$,
where $sV = V[1]$ is the graded vector space obtained by shifting the (homological) degree
of $V$ up by $1$.
The inclusion of $V$ into the two copies $V$ and $V'$ in $ \cyl(R) $ induces two maps of 
graded algebra $i : R \ra \cyl(R) $ and $i' : R \ra \cyl(R)$.

We say that a map $S : R \ra \cyl(R)$ of graded $k$-modules of degree $-1$ 
is a $(i,i')$-\emph{biderivation} if, for all homogeneous elements $a,b\in R$, we have
\[ 
S(ab) = S(a) \cdot i'(b) + (-1)^{|a|}i(a) \cdot S(b)
\]
It is easy to see that there exists a unique $(i,i')$-biderivation
$S : R \ra \cyl(R) $ such that $S(v) = sv \in sV \subset \cyl(R) $ for all $v \in V$.
This biderivation $S$  allows us to define a differential on $\cyl(R)$.
Indeed, there exists a unique derivation $d_{\cyl} : \cyl(R) \ra \cyl(R) $ 
of degree $-1$ such that
\begin{equation}   
\label{Baues_Lemaire_biderivation_definition}
\begin{split}
(1) \qquad &d_{\cyl} \circ i = i \circ d \\
(2) \qquad &d_{\cyl} \circ i' = i' \circ d \\
(3) \qquad &d_{\cyl} \circ S = i-i' - S \circ d
\end{split}
\end{equation}
Moreover, this derivation satisfies $\,d_{\cyl}^2 = 0\,$, which is easy to check on 
generators of $ \cyl(R) $. Hence, $\,d_{\cyl}\,$ turns the graded algebra $ \cyl(R) $
into a DG algebra $ \cyl(R) = (T(V \oplus V' \oplus sV), d_{\cyl})$.

Next, we define a map $\pi : \cyl(R) \ra R$ by sending the two copies of $V$  
in $R = T(V \oplus V' \oplus sV)$ identically onto $V \subset T(V) = R$ and $sV$ to zero. 
It is easy to check that this map commutes with differentials. Moreover, we have
$\, \pi \circ i = \id \,$ and
$\,\id - i\circ \pi = d_{\cyl} \circ S + S \circ d_{\cyl} \,$. Hence, $ \pi $ is 
a homotopy equivalence and therefore a quasi-isomorphism of DG algebras.

Thus, together with $\,i, i': R \ra \cyl(R) \,$, the map $ \pi $ fits in the
diagram $\,R \amalg R \xinto{(i,i')} \cyl(R) \stackrel{\pi}{\onto} R\,$,
which identifies $ \cyl(R) $ as a cylinder object of $R$.
We stress that this cylinder object is \emph{canonically} associated to the semi-free 
DG algebra $R$. We call it the \emph{Baues-Lemaire cylinder} on $R$.

The above construction can be naturally extended to semi-free DG categories, {\it
i.e.} DG categories whose underlying graded category is freely generated by a set of arrows.
In our present situation, 
the underlying graded category of the DG category $B$
is freely generated by the graded quiver \eqref{quiver_widetilde_Q_definition}.
Hence, the Baues-Lemaire construction of the cylinder on $R=T(V)$ can be carried over to 
$B = k\langle \tilde{Q} \rangle$.

Specifically, let $c\tilde{Q} = \tilde{Q} \amalg \tilde{Q}' \amalg (\tilde{Q}[1])$ be the graded quiver
\begin{equation}  \label{quiver_c_widetilde_Q_definition}
c\tilde{Q} = \qquad 
\vcenter{
\xymatrix{
\infinity   \ar@/^2.5pc/[rrrr]|{b}  \ar@/^1.5pc/[rrrr]|{a} \ar@/^0.5pc/[rrrr]|{a'}  
\ar@`{(-15,15),(-15,-15)}[]|{\xi}
\ar@`{(-20,20),(-20,-20)}[]|{\xi'}
\ar@`{(-25,25),(-25,-25)}[]|{\eta}
 & & & & 
0  \ar@/^2.5pc/[llll]|{b^*} \ar@/^1.5pc/[llll]|{a^*}  \ar@/^0.5pc/[llll]|{a'^*}
}
}
\end{equation}
which has three copies $\{a,a^*,\xi\}$, $\{a',a'^*,\xi'\}$
and $\{b,b^*,\eta\}$ of the generating arrows of $\tilde{Q}$,
with $\{b,b^*,\eta\}$ having homological degree shifted up by $1$. Thus,
\[
\begin{split}
\deg(a) &= \deg(a') = \deg(a^*) = \deg(a'^*) = 0\\
\deg(\xi) &= \deg(\xi') = 1 \\
\deg(b) &= \deg(b^*) = 1 \\
\deg(\eta) &= 2 \\
\end{split}
\]
Then, we define $\cyl(B) $ to be the graded $k$-category $\,\cyl(B) := k\langle c\tilde{Q} \rangle \,$,
with differential $ d = d_{\cyl} $ given by the Baues-Lemaire formulas \eqref{Baues_Lemaire_biderivation_definition}:
\begin{equation}
\begin{split}
 d(\xi) &= a^*a - (\mu-1)e_1 \\
 d(\xi') &= a'^* a' - (\mu-1)e_1 \\
 d(b) &= a-a' \\
 d(b^*) &= a^*-a'^* \\
 d(\eta) &= \xi - \xi' - b^*a' - a^*b
\end{split}
\end{equation}
For example, by Equation (3) in \ref{Baues_Lemaire_biderivation_definition}, we have
\[
\begin{split}
d(\eta) &= d(S(\xi)) = i(\xi) - i'(\xi) - S(d(\xi)) \\
&= \xi - \xi' - S(a^* a - (\mu-1)e_1) \\
&= \xi - \xi' - b^* a' - a^* b
\end{split}
\]

The following proposition implies that $\cyl(B)$ is indeed a cylinder object on $B$.
\bpp
The canonical map $\pi : \cyl(B) \ra B$ defined by 
\[
\begin{split}
\pi(a) &= \pi(a') = a, \qquad \pi(a^*) = \pi(a'^*) = a, \qquad \pi(\xi) = \pi(\xi') = \xi, \\
\pi(b) &= 0, \qquad \pi(b^*) = 0, \qquad \pi(\eta) = 0
\end{split}
\]
is a quasi-isomorphism.
\epp

Now, as explainedin the Introduction, the homotopy pushout \eqref{homotopy_closure_of_GMV_hocoeq_form}
can be computed as the ordinary pushout of the following diagram
\begin{equation}   \label{homotopy_closure_of_GMV_action_explicit_colimit}
 A^{(n)} \xla{(\Psi , \id^{(n)}) \circ (p^{(n)} , p^{(n)}) } 
 B^{(n)} \amalg B^{(n)} \xinto{(i^{(n)},i'^{(n)})} \cyl(B)^{(n)}
\end{equation}
A straightforward calculation shows that the result is the knot DG category presented in
Definition \ref{knot_DG_category_definition}.

\section{The knot category}   
\label{Section_cord_cat}
Let $ K $ be a knot, and let $ \cA_K $ be the knot DG category of $K$ presented in Definition 
\ref{knot_DG_category_definition}.
\bdf
We call the 0-th homology of $ \cA_K $ the \emph{knot $k$-category} of $K$ and 
denote it by $ A_K := H_0(\cA_K)$. This is a $k$-category whose isomorphism class is a knot invariant.
\edf

Let $\pi = \pi_1(\bR^3\!\setminus\! K)$ be the knot group of $K$.
Consider the group algebra $k[\pi]$ as a $k$-category with one object $0$.
Similarly, the ring $k$ can itself be considered a $k$-category with one object $1$,
which we denote by $ {\mathbf 1}_{\{1\}}$.
Let $k[\pi]^+ = k[\pi] \amalg {\mathbf 1}_{\{1\}}$ 
be the disjoint union of these two $k$-categories.
Thus, $k[\pi]^+ \in \Cat^{\{ 0,1 \}}_k$ is a $k$-category with object set $\{ 0,1 \}$.
Now, let $k[\pi]^+ \langle a,a^* \rangle$ be the free extension in $\Cat^E_k$
of $k[\pi]^+$ by the arrows $a,a^*$ where $a$ goes from the vertex $1$ to the vertex $0$,
while $a^*$ goes in the opposite direction, i.e., from $0$ to $1$.
We will denote this $k$-category schematically by
\[
k[\pi]^+ \langle a,a^* \rangle = 
\left[
\vcenter{
\xymatrix{
\bullet  \ar@/^/[rr]^{a}  
& & *+[F]{k[\pi]} \ar@/^/[ll]^{a^{*}}
}
}
\right]
\]
Then, we have the following description of the knot $k$-category 
in terms of the peripheral pair $(\pi, (m,l))$, where  
$m,l \in \pi$ are respectively a meridian and a longitude of the knot $K$.
\bthm  
\label{cord_category_main_theorem}
The knot $k$-category $A_K$ can be described as 
\[ A_K  \cong k[\pi]^+ \langle a,a^* \rangle / J \]
where $J$ is the ideal generated by the following elements
\begin{enumerate}
\item[$(1)$]\ $aa^* + e_0 - m$
\item[$(2)$]\ $a^*a + e_1 - \mu e_1$
\item[$(3)$]\ $\lambda a - l a \ ,\ \lambda a^* - a^* l$\ .
\end{enumerate}
\ethm

\brm
A peripheral pair $(m,l)$ is well-defined up to inner automorphisms of $\pi$.
Suppose that $(m',l')=(\gamma m \gamma^{-1}, \gamma l \gamma^{-1})$ is another such pair,
then letting $a' = \gamma a$ and $(a^*)' = a^* \gamma^{-1}$, we reduce the defining relations 
$(1)$-$(3)$ of the $k$-category $ A_K $ to the same form written
in terms of $a',(a^*)',l',m'$. Hence, up to isomorphism, 
this $k$-category is independent of the choice of the peripheral pair.
\erm

To prove Theorem~\ref{cord_category_main_theorem} we notice that the braid group $B_n$ acts on 
the elements $T_i \in A^{(n)}(0,0)$ the same way as it acts on the generators $x_i \in \bF_n$ 
in the Artin representation.
This implies that, after taking the categorical braid closure, there is a map $\phi$
from $k[\pi]$ to the endomorphism algebra $ A_K(0,0)$ of the knot $k$-category at $0$,
taking $x_i$ to $T_i$.
Define $\tilde{\phi} : k[\pi]^+ \langle a,a^* \rangle \ra A_K $ by
extending the map $\phi$, so that $a \mapsto a_1$, and $a^* \mapsto a_1^*$.
Then, one can show that if $ m = T_1 \in \pi$ and  $l \in \pi$ 
is the corresponding longitude, then the map $\tilde{\phi}$ sends the ideal $J$ defined in Theorem \ref{cord_category_main_theorem} 
to zero, and hence descends to a map from the quotient $k[\pi]^+ \langle a,a^* \rangle / J$ to $A_K$,
which can be shown to be an isomorphism. (See \cite{BEY} for details.)

Recall (see Theorem \ref{endomorphism_DGA_at_1_recovers_knot_DGA_theorem}) that the endomorphism DG algebra of the object $1$ in the knot DG category is quasi-isomorphic to the knot DGA.
Therefore, in particular, the endomorphism algebra of the object $1$ of the knot $k$-category 
recovers the $0$th homology of the knot DGA. Thus, Theorem \ref{cord_category_main_theorem} implies 

\bthm[\cite{Ng08}]
The $0$th homology of the knot DGA is isomorphic to the tensor algebra over $k$
freely generated by elements $[\gamma]$, where $\gamma \in \pi_1(\bR^3 \setminus L)$,
modulo the relations
\begin{enumerate}
 \item[$(1)$]\ $[e]=1-\mu$, where $e$ is the identity element;
 \item[$(2)$]\ $[\gamma_1 \gamma_2] - [\gamma_1 m \gamma_2] - [\gamma_1][\gamma_2]=0$ for $\gamma \in \pi_1(\bR^3 \setminus L)$;
 \item[$(3)$]\ $[\gamma l] = [l \gamma] = \lambda[\gamma]$ for $\gamma_1, \gamma_2 \in \pi_1(\bR^3 \setminus L)$.
\end{enumerate}
\ethm

\bpf
The endomorphism algebra of the object $1$ of the $k$-category $k[\pi]^+ \langle a,a^* \rangle$
is freely generated by the elements $[\gamma] := -a \gamma a^*$.
The ideal $J$ of Theorem \ref{cord_category_main_theorem} defines relations in this endomorphism algebra,
which are simply the three relations given in the theorem.
\epf

\brm
Theorem \ref{cord_category_main_theorem} also shows that the endomorphism algebra of the knot category 
at the vertex $0$  is given by 
$\, A_K(0,0) = k[\pi] /\langle (m-1)(m-\mu)\, ,\ (m-1)(l - \lambda) \rangle \,$.
%
\erm

\section{The fully noncommutative link DG category}
\label{FN}
Recall, in Section \ref{The_GMV_operator_section}, we have `simplified' the GMV $k$-category 
$\tilde{A}^{(n)} $ by performing the following two operations on the underlying quiver:
\begin{enumerate}
\item we have collapsed the vertices $1,\ldots,n$ to a single vertex $1$,
\item we have set the elements $ e_i + a_i^*a_i $ to be equal to 
a central element $\mu \in k^{\times}$.
\end{enumerate}

In this section, we will work with the original GMV category $\tilde{A}^{(n)} $ and
the associated braid action. We will show that the corresponding homotopy braid closure is related to the 
``fully noncommutative knot DGA'' introduced in \cite{EENS13a}, \cite{Ng14}.

Let $\tilde{A}^{(n)}$ be the $k$-category \eqref{widetilde_A}
with the GMV braid action defined as in \eqref{GMV}. Consider the elements
$ \mu_i = a_i^* a_i + e_i \in \tilde{A}^{(n)}(i,i)$, which are now no longer central.

Suppose we are given a braid $\beta \in B_n$ which closes to a link $L$ with $r$ component 
$L = L_1 \cup \ldots \cup L_r$.
For each $1\leq i \leq r$, let $S_i$ be the set of strands $S_i \subset \{1,\ldots,n\}$ 
that closes to the component $L_i$. Note that these
are precisely the orbits of the cyclic group generated by $\beta$ acting 
on the set $\{1,\ldots,n\}$ by permutations.

Now, for each $1 \leq i \leq r$,
identify all the vertices $j$ in $A^{(n)}$ that are in the set $S_i$ to a single vertex $i$,
and identify all the elements $\mu_j$, for $j\in S_i$, to a single element $\mu_i$.
Let $\overline{A^{(n)}}$ be the resulting $k$-category.
Then, the action map $\beta : A^{(n)} \ra A^{(n)}$ induces
$\overline{\beta} : \overline{A^{(n)}} \ra \overline{A^{(n)}}$.
We can use this induced braid action to define the fully noncommutative link DG category.

\bdf
\label{FNCDG}
The {\it fully noncommutative link DG category} of $L$ is the DG category $\tilde{\cA}_L $,
whose underlying graded $k$-cateogry is defined to be the quotient of
\begin{equation}
\xymatrix{
\xybox{
 *+[F]{ k[\lambda_1^{\pm 1}, \mu_1^{\pm 1}]  }
\ar@`{(-25,25),(-25,-25)}|{ \{ \eta_j \}_{j \in S_1} }
}
\ar@/^4pc/[rrrrd]|{ \{ b_j \}_{j \in S_1} }  
\ar@/^2.6pc/[rrrrd]|{ \{ a_j \}_{j \in S_1} }  
 & & & &  \\
\vdots & & & &  \bullet   
\ar@/_1.2pc/[llllu]|{ \{ b_j^* \}_{j \in S_1} }  
\ar@/^0.2pc/[llllu]|{ \{ a_j^* \}_{j \in S_1} } 
\ar@/^2.6pc/[lllld]|{ \{ b_j^* \}_{j \in S_r} }  
\ar@/^4pc/[lllld]|{ \{ a_j^* \}_{j \in S_r} } 
 \\
\xybox{
 *+[F]{ k[\lambda_r^{\pm 1}, \mu_r^{\pm 1}]  }
\ar@`{(-25,25),(-25,-25)}|{ \{ \eta_j \}_{j \in S_r} }
 }
\ar@/^0.2pc/[rrrru]|{ \{ b_j \}_{j \in S_r} }  
\ar@/_1.2pc/[rrrru]|{ \{ a_j \}_{j \in S_r} }  
}
\end{equation}
modulo the relations $e_i + a_j^* a_j = \mu_i$ for all $j \in S_i$, $1 \leq i \leq r$,
where the degrees of the generators are given by
\[
\deg(a_j) =  \deg(a_j^*) = 0, \qquad 
\deg(b_j) = \deg(b_j^*) = 1, \qquad
\deg(\eta_j) = 2
\]
To define the differential, choose a strand $j_i \in S_i$, one for each $1 \leq i \leq r$,
and for each $j \in S_i$, set
\begin{equation}
\begin{split}
 d(b_j) &= 
 \begin{cases} 
 \overline{\beta}(a_j ) \lambda_i^{-1} \mu_i^{-w_i} - a_j & (j = j_i) \\
 \overline{\beta}(a_j ) - a_j  & (j \neq j_i)
 \end{cases}
  \\
   d(b^*_j) &= 
 \begin{cases} 
  \lambda_i \mu_i^{w_i} \overline{\beta}(a^*_j ) - a_j & (j = j_i) \\
 \overline{\beta}(a^*_j ) - a_j  & (j \neq j_i)
 \end{cases}
 \\
 d(\eta_j) &= 
 \begin{cases} 
  - b^*_j a_j - \lambda_i \mu_i^{w_i} \overline{\beta}(a^*_j) b_j   & (j = j_i) \\
  - b^*_j a_j - \overline{\beta}(a^*_j) b_j   & (j \neq j_i)
 \end{cases}
\end{split}
\end{equation}
\edf
\bthm  
\label{fully_noncommutative_link_DG_category_main_theorem}
Let $R_0$ be the $k$-category with $r$ objects, given by the disjoint union of $k$-algebras
\[ R_0 = k[\lambda_1^{\pm 1}, \mu_1^{\pm 1}] \, \amalg \, \ldots \, \amalg k[\lambda_r^{\pm 1}, \mu_r^{\pm 1}]\]
Then, the quasi-isomorphism type of the pair $(R_0, \tilde{\cA})$ consisting of the $k$-category $R_0$, 
together with the canonical map from $R_0$ to the fully noncommutative link DG category $\tilde{\cA}$,
is a link invariant.
Moreover, if we collapse the objects $\{ 1, \ldots ,r \}$ to a single object $1$,
then the endomorphism DG algebra at this collapsed vertex coincides with the 
fully noncommutative knot DGA constructed in \cite{EENS13a}.
$($Here, we take the base commutative ring $k$ to be $\mathbb{Z}$.$)$
\ethm
The first part of the above theorem is proved by interpreting the fully noncommutative link DG category 
as a homotopy braid closure in a suitable model category. The second part follows from the first
by a direct calculation similar to the one in Section~\ref{Section_knot_DG_cat}
(see also the beginning of Section \ref{Generalizations_section}).
The identification of the fully noncommutative link DG category with a homotopy braid closure 
is completely parallel to the $\mu$-central case discussed above. 
The crucial difference, however, is that one should work in a different model category
(see \cite{BEY} for details).

The above theorem identifies the \emph{quasi-isomorphism} type of the pair $(R_0, \tilde{\cA})$; however,
if we are only interested in the underlying \emph{quasi-equivalence} type, then we have the following

\bthm
\label{TTT}
The quasi-equivalence type of the link DG category $\tilde{\cA}$ 
is given by the (writhe-adjusted) homotopy closure of the braid $\beta \in B_n$ 
with respect to the GMV operator,
taken in the category $\dgcat^*_k$ of pointed DG categories
with  model structure defined in \cite{Tab05}. 
\ethm

Notice that, in this theorem, no coloring is needed.
The extra parameters $\lambda_i$ are formed in the process of taking the homotopy braid closure.
This is not ``visible'' if we, like in the $\mu$-central case, 
work with a more rigid model structure, where the weak equivalences are quasi-isomorphisms
({\it cf.} also Remark \ref{catcl_of_GMV_is_not_cord_cat}).
``Writhe-adjusting'' is also not necessary in Theorem~\ref{TTT}, as it only changes the parameter 
$\,\lambda_i \mapsto \lambda_i \,\mu_i^{w_i}\,$ in $R_0$.

\bdf
The {\it fully noncommutative link $k$-category} of a link $L$ is defined to be 
$ \tilde{A}_L := H_0(\tilde{\cA}_L) $, the $0$th homology of the fully 
noncommutative link DG category of $L$.
\edf
The $k$-category $ \tilde{A}_L $
can be expressed in terms of the link group, together with meridians and longitudes
chosen in each link component. To be precise, let $M = \bR^3\! \setminus\! L $ be the link complement.
For $1 \leq i \leq r$, let $\partial_i M \subset M$ denote
the torus boundary of $M$ corresponding to the link component $L_i$.
Choose basepoints $p_i \in \partial_i M$, and $p_0 \in M$.
Then, there are canonical meridian and longitude elements $\mu_i, \lambda_i \in \pi_1(\partial_i M, p_i)$,
which identifies the group algebra $k[\pi_1(\partial_i M, p_i)]$ as $k[ \lambda_i^{\pm 1}, \mu_i^{\pm 1}]$.
By choosing a path $a_i$ in $M$ from $p_i$ to $p_0$, one can define a map 
$\phi_i : \pi_1(\partial_i M, p_i) \ra \pi_1(M, p_0)$.
Let $m_i$ and $l_i$ be the images of $\mu_i$ and $\lambda_i$ under $\phi_i$, respectively.
Then, we have the following description of the fully noncommutative link category.

\bthm  
\label{noncommutative_cord_category_peripheral_description}
The fully noncommutative link $k$-category $ \tilde{A}_L$ is 
the quotient of the $k$-category
\begin{equation}
\xymatrix{
 *+[F]{ k[\lambda_1^{\pm 1}, \mu_1^{\pm 1}]  }
\ar@/^2pc/[rrrrd]|{  a_1  }  
 & & & &  \\
\vdots & & & &  
*+[F]{ k[ \pi_1(M,p_0) ]  }
\ar@/^1pc/[llllu]|{ a_1^* } 
\ar@/^2pc/[lllld]|{ a_r^* } 
 \\
 *+[F]{ k[\lambda_r^{\pm 1}, \mu_r^{\pm 1}]  } 
\ar@/^1pc/[rrrru]|{  a_r  }  
}
\end{equation}
modulo the ideal of relations
\begin{enumerate}
\item[$(1)$]\ $a_ia_i^* + e_0 - m_i$
\item[$(2)$]\ $a_i^*a_i + e_i - \mu_i$
\item[$(3)$]\ $a_i \lambda_i - l_i a_i \ ,\ \lambda_i a_i^* - a_i^* l_i$
\end{enumerate}
\ethm

\brm  
\label{catcl_of_GMV_is_not_cord_cat}
While the fully noncommutative link DG category is the homotopy braid closure of the GMV action,
it is \emph{not} true that its $0$-th homology, {\it i.e.} 
the fully noncommutative link $k$-category, is the categorical braid closure of the GMV action.
The categorical braid closure can be obtained as a specialization of the fully noncommutative link $k$-category 
when all parameters $\lambda_i$ are set to be $1$.
This discrepancy is due to the fact that, in Tabuada's model structure on $\dgcat_k$,
the weak equivalences are quasi-equivalences, which, by definition, induce 
equivalences (not isomorphisms) of $k$-categories at the level of $0$-th homology.
Then, the homotopy colimits of diagrams in $\dgcat_k$ with respect to Tabuada's model structure 
induce, at the level of $0$th homology, not strict colimits, but rather $2$-colimits, which can be viewed, in part, 
as homotopy colimits. Thus, the fully noncommutative link $k$-category is already a homotopy braid closure,
rather than a strict categorical braid closure.
\erm

As mentioned in the introduction, the fully noncommutative link $k$-category 
is closely related to perverse sheaves. To be precise, let $\cS$ be the stratification on 
$\bR^3$ with two strata $(L,\bR^3\!\setminus\! L)$, where $L$ is a link in $ \bR^3 $.
Following the degree conventions of \cite{KS94}, we let $p$ be the perversity of $ \cS $ given by $p(1)=0$ 
and $p(3) = -1$. (The values at other integers do not matter.)
Then, we have
\bthm  
\label{relation_with_perverse_sheaves}
Suppose that $k$ is a field.
The category $\Perv^p(\bR^3,L)$ of $p$-perverse sheaves of $k$-vector spaces on $\bR^3$ 
constructible with respect to the 
stratification $\cS$ is equivalent to the category of finite-dimensional left modules over the 
fully noncommutative link category $\tilde{A}_L $.
\ethm

\begin{proof}[Sketch of proof]
Suppose that a braid $\beta \in B_n$ is placed in the region $\{x < 0\}$,
and closes to the link $L$ by letting the two ends of the braid pass through the hyperplane $\{x=0\}$
and close in the region $\{x>0\}$, as in the following diagram.
\[ 
\includegraphics[scale=0.3]{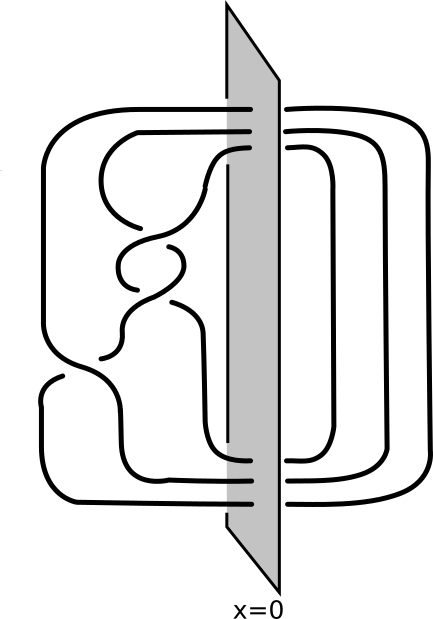} 
\]
Let $U,V$ be open subsets of $\bR^3$ defined by $U = \{ x<\varepsilon \}$ and $V = \{ x > - \varepsilon \}$
for some small $\varepsilon >0$.
Then, both the pairs $(U,U\cap L)$ and $(V,V\cap L)$ are diffeomorphic to the pair 
$(\mathring{D} \times \mathring{I}, \{p_1,\ldots,p_n\} \times \mathring{I} )$,
where $\mathring{D}$ denotes the interior of the disk and $\mathring{I}$ denotes the open interval $(0,1)$.
The pair $(U\cap V, U\cap V \cap L)$ is diffeomorphic to the pair 
$(\mathring{D} \times \mathring{I}, \{p_1,\ldots,p_n, p'_1, \ldots, p'_n\} \times \mathring{I} )$.
Therefore, the category $\Perv^p(U,U\cap L)$ can be identified with the category 
$\Perv(D, \{p_1, \ldots ,p_n\})$ with middle perversity,
which, under a suitable choice of `cuts', is  equivalent to the category 
$\Mod(\tilde{A}^{(n)})$ of finite-dimensional modules over the $k$-category $\tilde{A}^{(n)}$.

Similar statements are true for the pairs $(V, V\cap L)$ and $(U\cap V, U\cap V \cap L)$.
One can show then that the following diagram of restriction functors
\begin{equation}  \label{perverse_sheaves_limit_diagram}
\Perv^p(U, U\cap L ) \ra  \Perv^p(U \cap V, U\cap V \cap L ) \la \Perv^p(V, V\cap L )
\end{equation} 
is equivalent to the following diagram of functors
\begin{equation}  \label{perverse_modules_limit_diagram}
\Mod(A^{(n)}) \xra{(\beta^* , \id)} \Mod(A^{(2n)}) \xla{(\id,\id)} \Mod(A^{(n)})  
\end{equation}

Since perverse sheaves form a stack (see \cite[Propositions 10.2.7 and 10.2.9]{KS94}), 
the category $\Perv^p(\bR^3,L)$ is equivalent to 
the $2$-limit of the diagram \eqref{perverse_sheaves_limit_diagram}, 
and hence of the diagram \eqref{perverse_modules_limit_diagram}.
This implies the desired result. For details, see \cite{BEY}.
\end{proof}
When combined with Theorem~\ref{noncommutative_cord_category_peripheral_description},
Theorem~\ref{relation_with_perverse_sheaves} gives a description of the category $\Perv(\bR^3,L)$ of perverse sheaves 
in terms of linear algebra data, similar in spirit to the original description of the category 
$ \Perv(D,\{p_1, p_2, \ldots, p_n\})$ given in \cite{GMV96}.

\section{Generalizations and further questions}  
\label{Generalizations_section}
In the GMV braid action, the group $B_n$ acts on the generators 
$T_i = a_i a_i^* + e_0$ via the Artin representation \eqref{artinact}.
Thus, regarding the free group $ {\mathbb F}_n $ as a category with a single object, 
we can regard the GMV action as an extension of the Artin action.
In \cite{Wad92}, Wada constructed several examples of braid group actions on 
$ {\mathbb F}_n $ generalizing the classical Artin representation. Like the
Artin representation, Wada's braid group actions are local and homogeneous, {\it i.e.}
generated by a single cocartensian Yang-Baxter operator on $ {\mathbb F}_2 = {\mathbb F}_1 \amalg {\mathbb F}_1 $.
It is therefore natural to ask whether they admit extensions similar to the GMV extension.

Consider, for example, the following cocartesian Yang-Baxter operator constructed in \cite{Wad92}:
\[ 
\sigma : \bF_2 \ra \bF_2 \qquad x_1 \mapsto x_1^N x_2 x_1^{-N}, \qquad x_2 \mapsto x_1\ ,
\]
where $ N $ is an arbitrary (fixed) integer. 

This action does admit an extension similar to the GMV action. Indeed,
\begin{equation} 
\label{GMV_N_version_Bn_action_equation}
\sigma_k \, : \,  \begin{cases} 
     a_i \, \mapsto \, a_{i} & \mbox{$(i \neq k, k+1)$} \\
     a_{k} \, \mapsto \, T_{k}^N a_{k+1} &  \\
     a_{k+1} \, \mapsto \, a_{k} & \\ 
      a_i^* \, \mapsto \, a_{i}^* & \mbox{$(i \neq k, k+1)$} \\
     a_{k}^* \, \mapsto \,  a_{k+1}^*  T_{k}^{-N}&  \\
     a_{k+1}^* \, \mapsto \, a_{k}^* &
   \end{cases}  \end{equation}
Note that, for $ N = 1 $, this is the original GMV action \eqref{GMV}.
Moreover, using a result of \cite{CP05}, one can show that the actions \eqref{GMV_N_version_Bn_action_equation}
are non-equivalent to each other for different $N$'s; thus, for $N \not= 1 $, \eqref{GMV_N_version_Bn_action_equation} 
is a genuine generalization of the GMV action.

The elements $T_i^{\pm N}$ can be written in an alternative form 
involving the conjugate elements $\mu_i = e_i + a_i^* a_i \in \tilde{A}^{(n)}(i,i)$.
(We recall that $ \mu_i $ are no longer central elements in $\tilde{A}^{(n)}(i,i)$.).
Indeed, by an induction, one can show that
\[ T_i^{N} = e_0 + a_i \,[N]_{\mu_i}\, a_i^*\qquad \text{for all}\ N\in \bZ\]
where $\,[N]_{\mu} \in k\,$ are the ``quantum integers'' defined by
\[ [N]_{\mu_i} = \frac{\mu_i^N - 1}{\mu_i - 1}
:= \begin{cases}
e_i + \mu_i + \mu_i^2 + \ldots + \mu_i^{N-1} & \text{if }N>0 \\
0 & \text{if }N=0 \\
-\mu_i^{-1} - \mu_i^{-2} - \ldots - \mu_i^{N} & \text{if }N<0
\end{cases}  \]
As in Section~\ref{The_GMV_operator_section}, we set $\,A_{ij} = -a_i^* a_j\,$ for all 
$i,j$. Then, we have the following formulas defining the braid group action 
on the restriction of the $k$-category $\tilde{A}$ to the vertices $\{1, \ldots ,r\}$:
\[
\sigma_k \, :  \, \begin{cases}
A_{ki} \, \mapsto \, A_{k+1,i} -  A_{k+1, k} \, [-N]_{\mu_k} \,  A_{ki} & i \neq k, k+1\\
A_{ik} \, \mapsto \, A_{i, k+1} -  A_{ik} \, [N]_{\mu_k} \, A_{k,k+1} & i \neq k, k+1\\
A_{k+1,i} \, \mapsto \, A_{ki}  & i \neq k, k+1 \\
A_{i,k+1} \, \mapsto \, A_{ik} &  i \neq k, k+1\\
A_{k, k+1}  \, \mapsto \,  A_{k+1, k} \, \mu_k^{-N}   &  \\
A_{k+1, k}  \, \mapsto \, \mu_k^N \, A_{k, k+1} & \\
A_{ij} \, \mapsto \, A_{ij}  &  i,j \neq k, k+1 \\
\mu_k \, \mapsto \, \mu_{k+1} \\
\mu_{k+1} \, \mapsto \, \mu_k \\
\mu_i \, \mapsto \, \mu_i & i \neq k, k+1
 \end{cases}
\]

Now, for $i \neq j$, define
%
\begin{equation} \label{Aij_aij_general_N_relation_substitution}
a_{ij} = 
\begin{cases}
A_{ij}\, [N]_{\mu_j} , &i<j \\
A_{ij}\, [-N]_{\mu_j}, &i>j
\end{cases}
\end{equation}
Then, the above action becomes
\[
\sigma_k \, :  \, \begin{cases}
a_{ki} \, \mapsto \, a_{k+1,i} -  a_{k+1, k}  \,  a_{ki} & i \neq k, k+1\\
a_{ik} \, \mapsto \, a_{i, k+1} -  a_{ik} \, a_{k,k+1} & i < k\\
a_{ik} \, \mapsto \, a_{i, k+1} -  a_{ik} \, \mu_k^{N} \, a_{k,k+1} \, \mu_{k+1}^{-N} & i > k+1\\
a_{k+1,i} \, \mapsto \, a_{ki}  & i \neq k, k+1 \\
a_{i,k+1} \, \mapsto \, a_{ik} &  i \neq k, k+1\\
a_{k, k+1}  \, \mapsto \,  - a_{k+1, k}   &  \\
a_{k+1, k}  \, \mapsto \, -\mu_k^N \, a_{k, k+1} \, \mu_{k+1}^{-N} & \\
a_{ij} \, \mapsto \, a_{ij}  &  i,j \neq k, k+1 \\
\mu_k \, \mapsto \, \mu_{k+1} \\
\mu_{k+1} \, \mapsto \, \mu_k \\
\mu_i \, \mapsto \, \mu_i & i \neq k, k+1
 \end{cases}
\]
For $N=1$, this coincides with the `fully noncommutative' action defined
in \cite{EENS13a} (see also \cite[Appendix]{Ng14}).

In a different direction, one can also construct a large family of 
GMV-type braid actions by extending the family of generalized Artin actions found in \cite{CP05}.
Specifically, let $B \in \Alg_k$ be an associative algebra over $k$, and let $x,y \in B^{\times}$ 
be a pair of \emph{invertible} and \emph{commuting} elements.
Let $\hatA^{(n)}$ be the $k$-category given by
\[
\hatA^{(n)} = 
\left[
\vcenter{
\xymatrix{
\bullet \ar@/^0.5pc/[rrd]^{a_1} & & \cdots & & \bullet \ar@/^0.5pc/[lld]^{a_n} \\
& & *+[F]{B\, * \stackrel{n}{\ldots} * \, B} \ar@/^0.5pc/[llu]^{a_1^*}  \ar@/^0.5pc/[rru]^{a_n^*}
}
}
\right]
\]
which can be interpreted as an $n$-fold coproduct
in the category $\Cat^*_k$ of (small) pointed $k$-categories.
Then, one can check by a direct calculation that the following assignments
define a braid group action on $\hatA^{(n)}$.
\begin{equation}  \label{Crisp_Paris_action}
\sigma_k \, : \,  \begin{cases} 
 a_{k} \, \mapsto \, x_k a_{k+1} &  \\
 a_{k}^* \, \mapsto \,  a_{k+1}^*  x_{k}^{-1}&  \\
 b_{k} \, \mapsto \,  x_k b_{k+1} x_{k}^{-1}&  \\ 
 a_{k+1} \, \mapsto \, y_k a_{k} & \\ 
 a_{k+1}^* \, \mapsto \, a_{k}^* y_k^{-1} & \\
 b_{k+1} \, \mapsto \, y_k b_k y_k^{-1} & \\ 
 a_i \, \mapsto \, a_{i} & \mbox{ $(i \neq k, k+1)$} \\
 a_i^* \, \mapsto \, a_{i}^* & \mbox{$(i \neq k, k+1)$} \\
 b_i \, \mapsto \, b_{i} & \mbox{$(i \neq k, k+1)$} 
   \end{cases}
\end{equation}
where $b_i$ is the element $b\in B$ put in the $i$-th copy of $B$ in 
$B^{(n)} := B\, * \stackrel{n}{\ldots} * \,B\,$.
Notice that, when $B = k[H]$ is the group algebra of a group $H$,
and when $x = h \in H$ and $y = h^{-1}$,
the braid action on $H^{(n)} \subset k[H]^{(n)}$ at the vertex $0$ 
coincides with the action defined in \cite{CP05}.

Consider any ideal $I \subset \hatA(0,0) = B\langle aa^* \rangle$,
and let $\hatA / I$ be the $k$-category obtained by quotienting $\hatA$ by the 
ideal generated by $I$.
Then, for any element $f \in I$, we have
\[ \sigma_k(f_k) = x_k f_{k+1} x_k^{-1}\qquad \sigma_k(f_{k+1}) = y_k f_k y_k^{-1} 
\qquad \sigma_k(f_i = f_i \text{ if }i \neq k,k+1)\]
Therefore, the cocartesian Yang-Baxter operator corresponding to the above braid group action
descends to the quotient $\hatA / I$.
If we take $B = k[T^{\pm}]$, $x = T$ and $y=1$,
and consider the ideal $I$ generated by the element $a a^* + 1 - T$,
then the resulting quotient $\hatA / I$, together with its corresponding cocartesian Yang-Baxter operator $\sigma$,
is equivalent to the $k$-category 
$\tilde{A}$, together with the GMV operator, constructed in Section \ref{The_GMV_operator_section}.

%
\bthm
The braid group actions \eqref{GMV_N_version_Bn_action_equation} and 
\eqref{Crisp_Paris_action} are generated by Reidemeister operators 
in the category $\dgcat^*_k$ of pointed DG categories on objects $ \tilde{A} $ and $ \hatA/I $, 
respectively. These objects are pseudoflat with respect to Tabuada's model structure on $\dgcat^*_k$.
Thus, the homotopy braid closure with respect to these operators gives 
link invariants which generalize the fully noncommutative link DG category $\tilde{\mathscr{A}}_L$.
\ethm
We conclude the paper with a few  questions and remarks.

\vspace{1ex}

{\bf 1.} {\it Constructible sheaves and contact homology.}
Recently, some interesting work has been done on the geometric side of contact homology
relating it to constructible sheaves (see \cite{ENS, NRSSZ, STZ, She, ST}). It would be
interesting to understand our Theorem \ref{relation_with_perverse_sheaves} in this 
geometric context and more generally, to clarify the meaning of our construction from 
Floer-theoretic and constuctible sheaves point of view. 

In more detail, the relation between Legendrian contact homology and constructible sheaves is based on 
a theorem of Nadler and Zaslow \cite{NZ09} (see also \cite{Nad09}) that, 
for any real analytic manifold $M$, establishes an equivalence between 
the derived category $ D_c(M) $ of constructible sheaves on $M$ and the derived 
Fukaya category $ D\Fuk(T^*M)$ of the cotangent bundle $T^*M$ of $M$. This equivalence of 
triangulated categories is induced by a quasi-equivalence of $A_{\infty}$-categories
$ \mu:\, \Sh_c(M) \ra \Tw\Fuk(T^*M) \,$, where $ \Sh_c(M) $ is a DG category defined as the DG
quotient of the (naive) DG category of constructible sheaves on $M$ modulo acyclic complexes 
and $\Tw\Fuk(T^*M)$ is the $A_{\infty}$-category  of twisted complexes in the Fukaya category 
$\Fuk(T^*M)$. The functor $\mu $ can be viewed as a categorification of the classical 
characteristic cycle construction and is called the \emph{microlocalization functor}. 

Now, for any conical Lagrangian submanifold $\tilde{\Lambda} \subseteq T^*M$, the restriction of
the microlocalization functor to the subcategory $\Sh_c(M)_{\tilde{\Lambda}} \subseteq \Sh_c(M)$ 
of constructible sheaves with singular support in $\tilde{\Lambda}$ gives a quasi-equivalence 
$ \mu: \Sh_c(M)_{\tilde{\Lambda}} \stackrel{\sim}{\to} \Tw\Fuk(T^*M)_{\tilde{\Lambda}}$ onto the full subcategory
$\Tw\Fuk(T^*M)_{\tilde{\Lambda}} $ of the twisted Fukaya category consisting of Lagrangians whose 
boundary at infinity lies in the boundary of $\tilde{\Lambda}$. Such a submanifold $ \tilde{\Lambda} $ is 
determined by its intersection $\, \Lambda := \tilde{\Lambda}\, \cap \,ST^*M \,$ with the unit cotangent bundle 
of $M$; the bundle $ ST^* M $ has a natural contact structure, and $ \Lambda $ is a Legendrian submanifold of $ST^*M $.
It turns out that the Legendrian contact homology (LHC) of the pair $(ST^*M, \Lambda)$ 
is related to the Fukaya category $ \Tw\Fuk(T^*M)_{\tilde{\Lambda}} $ 
and hence, via the microlocalization functor, to the sheaf category $ \Sh_c(M)_{\tilde{\Lambda}}$.
More precisely, it is expected that the complexes of constructible sheaves in $\Sh_c(M)_{\tilde{\Lambda}}$ 
determine augmentations of the Legendrian DGA of $(ST^*M, \Lambda)$ via a geometric symplectic filling
construction.

In the case of one-dimensional Legendrians, this relation has been worked out in detail in \cite{NRSSZ, STZ}.
Specifically, if $\,M = \bR^2\,$, then $ST^*\bR^2 \cong \bR^2 \times S^1$ contains an open contact submanifold 
$\bR^3 \subset \bR^2 \times S^1$. Hence, any Legendrian link $L \subset \bR^3$ can be considered 
as a Legendrian submanifold in $ ST^*\bR^2  $. In \cite{NRSSZ}, for a Legendrian link $L \subset \bR^3$, 
the authors construct a (unital) $A_{\infty}$-category $\Aug_{+}(L)$, whose objects are augmentations of the Chekanov-Eliashberg DG algebra of $L$, and show that there is an $A_{\infty}$-equivalence $\,\Aug_{+}(L) \simeq \cC_1(L)\,$, 
where $\cC_1(L)$ is the full subcategory of $ \Sh_c(\bR^2)_{\tilde{L}}$ consisting of sheaves of 
`microlocal rank one along the link $L$'.

A possible extension of this equivalence to higher dimensions (specifically, 
to the case of knot contact homology and knot DGA in $ \bR^3 $) has been recently proposed by
V.~Shende et. al. (see, e.g., \cite[Section 4]{She}, \cite[Section 6.6]{ENS}, \cite[Section 6.5]{ST}).
In this case, $M = \bR^3$ and the Legendrian $ \Lambda \subset ST^*M $ is given by the unit conormal bundle 
$\Lambda_L := ST^*_L \bR^3 $ associated to a link $L \subset \bR^3$. It is interesting that the support 
condition defining the subcategory $\Sh_c(\bR^3)_{\tilde{\Lambda}} \subset \Sh_c(\bR^3) $ coincides with the 
constructibility condition in our Theorem \ref{relation_with_perverse_sheaves}, and some 
geometric arguments suggest that there is a relation between this sheaf category and knot contact 
homology (see  \cite[Section 6.6]{ENS}). Whether this geometric relation can be used to prove the result
of Theorem~\ref{relation_with_perverse_sheaves} is not clear to us at the moment: {\it a priori}, the 
equivalence of categories in Theorem~\ref{relation_with_perverse_sheaves} originates from a different 
direction. In fact, there are three approaches to knot contact homology:
\begin{enumerate}
\item combinatorial knot contact homology,
\item Legendrian contact homology of the pair $\Lambda_L \subset ST^* \bR^3$,
\item constructible sheaves on $\bR^3$ with singular support in $ \tilde{\Lambda}_L $.
\end{enumerate}
The papers \cite{EENS13a, EENS13b} establish an equivalence between $(1)$ and $(2)$
by identifying the generators of the combinatorial knot DGA 
with Reeb cords and defining the differentials in terms of pseudoholomorphic curves. 
The geometric approach of \cite{ENS, She, ST} relates $(2)$ and $(3)$ 
via the geometry of symplectic fillings. Our result, Theorem \ref{relation_with_perverse_sheaves}, establishes 
the relation between $(1)$ and $(3)$ by appealing to the classical description of perverse sheaves on the disk
in terms of nearby and vanishing cycle functors \cite{GMV96} and using an algebraic `gluing' construction 
(homotopy braid closure). It would be interesting to see whether these approaches actually `agree'; 
in particular, can one prove Theorem~\ref{relation_with_perverse_sheaves} using the approach of \cite{ENS, She, ST}?

{\bf 2.} {\it Categorification of the  link DG category.}
There seems to be a natural way to categorify the DG category $ \tilde{\cA} $, using
the notion of `perverse schobers' introduced in \cite{KS1} (see also \cite{KS2}).
First, one can construct a (higher) category $\mathscr{C}$  of
$(\infty,2)$-categories that includes the category 
$\dgcat_k$ as an object (see \cite{Tam07, Fao}).
In $\mathscr{C}$, one can find an object $\mathscr{A}^{(n)}$ 
such that the category of $2$-representations of $\mathscr{A}^{(n)}$, 
{\it i.e.} an appropriately defined internal hom $\underline{\Hom}_{\mathscr{C}}( \mathscr{A}^{(n)} , \dgcat_k )$,
is equivalent to the (higher) category of perverse schobers on the disk with $n$ marked points.
Then, there should exist a $B_n$-action on $\mathscr{A}^{(n)}$ for all $ n \ge 1 $ that would allow us to
take the homotopy braid closure.
The result should be an object in $\mathscr{C}$ ({\it i.e.}, an $(\infty,2)$-category $\mathbb{A}$), whose 
category $\underline{\Hom}_{\mathscr{C}}( \mathbb{A} , \dgcat_k )$ of $2$-representations
is equivalent to a category of `perverse schobers on $\bR^3$ singular along a link'.

{\bf 3.}
{\it Yang-Baxter operators related to coherent sheaves.}
Many interesting examples of braid group actions related to coherent sheaves have been
constructed in the literature (see, e.g., \cite{ST01}, \cite{Rou06},\cite{AL1} and references therein).
It would be interesting to look at these examples in relation to the examples studied in 
the present paper and clarify the relations between the corresponding link invariants.



\end{document}